\documentclass[reqno,11pt]{amsart}

\usepackage{amsthm,amsmath,amssymb}
\usepackage{mathrsfs,amsfonts,dsfont,functan,extarrows,mathtools}

\usepackage{hyperref}
\usepackage{xcolor}

\newif\ifarxiv
\arxivtrue        

\ifarxiv
  \hypersetup{
      colorlinks=true,
      linkcolor=[rgb]{0.1,0.2,0.6},   
      citecolor=[rgb]{0.1,0.2,0.6},   
      urlcolor=[rgb]{0.1,0.2,0.6},    
      anchorcolor=[rgb]{0.1,0.2,0.6},
      pdfborder={0 0 0},              
      pdftitle={First eigenvalue estimates on complete balanced Hermitian manifolds},
      pdfauthor={Liangdi Zhang},
      pdfsubject={53C55; 58C40},
      pdfkeywords={balanced Hermitian manifold; Strominger--Bismut connection; Laplace--de~Rham operator; holomorphic Ricci curvature; holomorphic sectional curvature; eigenvalue estimate.}
  }
\else
  \hypersetup{
      colorlinks=false,               
      pdfborder={0 0 0},              
      pdftitle={First eigenvalue estimates on balanced Hermitian manifolds},
      pdfauthor={Liangdi Zhang}
  }
\fi

\usepackage{marginnote}
\usepackage{xcolor}

\makeatletter

\newcommand{\Rmnum}[1]{\expandafter\@slowromancap\romannumeral #1@}
\makeatother

\newtheorem{theorem}{Theorem}[section]
\newtheorem{proposition}[theorem]{Proposition}

\newtheorem{corollary}[theorem]{Corollary}

\newtheorem{lemma}[theorem]{Lemma}
\newtheorem{example}[theorem]{Example}
\newtheorem{conjecture}[theorem]{Conjecture}

\numberwithin{equation}{section}
\allowdisplaybreaks

\begin{document}

\title[]{First eigenvalue estimates on complete balanced Hermitian manifolds}

\author[]{Liangdi Zhang}
\address[Liangdi Zhang]{\newline Mathematical Science Research Center, Chongqing University of Technology, Chongqing 400054, China}
\email{ldzhang91@163.com}



\begin{abstract}

We establish lower bounds for the first positive eigenvalue of the
Laplace--de Rham operator on complete balanced Hermitian manifolds in
terms of curvature of the Strominger--Bismut connection. Under a positive
lower bound for its holomorphic Ricci curvature, we prove a
Lichnerowicz--Obata type estimate and characterize the equality case in the
K\"ahler setting. We also derive Li--Yau and Zhong--Yang type estimates from
lower bounds on the same holomorphic Ricci curvature, including estimates
under weaker assumptions only along a first eigendirection in the compact
case. Finally, under a positive lower bound for the holomorphic sectional
curvature of the Strominger--Bismut connection and a torsion-commutator
condition along a first eigendirection, we obtain a lower bound for the
first eigenvalue. In the compact case, the commutator condition follows from
vanishing of the Strominger--Bismut torsion in that eigendirection. These
results extend several classical K\"ahler and Riemannian spectral estimates
to the balanced non-K\"ahler setting.

\vspace*{5pt}


\end{abstract}

\maketitle

\tableofcontents

\section{Introduction}

Let $(M,\omega)$ be a complete Hermitian manifold of complex dimension $n$ with the fundamental $2$-form
\[
\omega = \sqrt{-1}\, h_{i\bar{j}} dz^i \wedge d\bar{z}^j.
\]
The Hermitian metric $h$ is called balanced if
\[
d\omega^{n-1} = 0,
\]
which is equivalent to $d^*\omega = 0$. The balanced condition is natural in this context, as it guarantees favorable analytical properties for the Laplace--de~Rham operator and aligns with the torsion structure of the Strominger--Bismut connection. In particular, any balanced Hermitian manifold is automatically K\"ahler in complex dimensions one and two.

Throughout this paper, let $\lambda_1>0$ denote the first (positive) eigenvalue of the Laplace--de~Rham operator $\Delta_d = d^* d + d d^*$
on a compact balanced Hermitian manifold $(M,\omega)$, i.e., there exists $u \in C^\infty(M,\mathbb{R})$ such that
\begin{equation}\label{1}
\Delta_d u = \lambda_1 u.
\end{equation}

Denote the diameter of $(M,\omega)$ by $D$ and the dual $(1,0)$-vector field of $\partial u = \frac{\partial u}{\partial z^k} dz^k \in \Gamma(M,T^{*1,0}M)$ by
\[
U = (\partial u)^\sharp = h^{i\bar{k}} \frac{\partial u}{\partial \bar{z}^k} \frac{\partial}{\partial z^i} \in \Gamma(M,T^{1,0}M).
\]

Classical lower bounds for the first eigenvalue in the Riemannian setting, including the works of Lichnerowicz--Obata \cite{Lz,O}, Cheng \cite{Cheng75}, Li--Yau \cite{LY80}, and Zhong--Yang \cite{ZY84}, provide fundamental tools in geometric analysis. Extensions have been developed for weighted manifolds \cite{BCP14,BH24,DMWX21}, manifolds with Bakry--\'{E}mery Ricci curvature bounds \cite{BQ00,xLi05,Wu10,Wu13,WSZ24}, manifolds with boundary \cite{CZL12,CY05,CY09,Gr02,LW21,yLY25,So02}, and pseudo-Hermitian manifolds \cite{CC07,CC09,Gre85,LL04,LT11,MV23}.

In K\"ahler geometry, several first eigenvalue estimates have been established, reflecting the interplay between the complex structure and the Riemannian metric. For a complete K\"ahler manifold $(M,g)$ of complex dimension $n$, Futaki \cite{Fut88} showed that $\lambda_1 \ge 2(n+1)$
when the holomorphic Ricci curvature satisfies $Ric(g) \ge (n+1)g$. Some early foundational results can be found, for instance, in \cite{LW05,Liu14,Mun09,Mun10,TY12}. More recently, Chu--Wang--Zhang \cite{CWZ25} proved a rigidity result: if in addition $\lambda_1 = 2(n+1)$ and the holomorphic bisectional curvature is positive, then $(M,g)$ is isometric to $(\mathbb{C}\mathbb{P}^n, \omega_{FS})$. In 2025, Wang--Yang \cite{WY25+} obtained a lower bound
\[
\lambda_1 \ge \frac{320(n-1)+576}{81(n-1)+144}
\]
under the condition $\mathrm{HSC} \ge 2$, using a Bochner--Kodaira type identity specifically developed for holomorphic sectional curvature. These results motivate our study of the analogous problem in the Hermitian, non-K\"ahler setting.

In Hermitian geometry, the choice of connection plays a crucial role. Among canonical Hermitian connections, the Strominger--Bismut connection ${^{SB}}\nabla$ arises naturally in complex geometry \cite{Bis89} and theoretical physics \cite{Str86}. Over the past decade, significant advances have been made in understanding the Strominger--Bismut connection; see, e.g., \cite{AV22,AOUV22,ABLS25,FTV22,NZ23,WY25,WYZ20,YZ18a,YZ18b,Yang25a+,YZZ23,Ye25,Zhang25+,ZZ19,ZZ23}.

Precise definitions and formulas for the Strominger--Bismut curvature are recalled in Section~\ref{sec2}. The holomorphic Ricci curvature of the Strominger--Bismut connection serves as a natural extension of classical Riemannian Ricci curvature to Hermitian geometry, while its holomorphic sectional curvature provides a conditional Hermitian extension of the Wang--Yang \cite{WY25+} estimate under a first-eigendirection torsion condition.

To the best of our knowledge, no Lichnerowicz--Obata type, Li--Yau type, or Zhong--Yang type eigenvalue estimates have been established on complete balanced Hermitian manifolds under curvature conditions associated with the Strominger--Bismut connection. In this paper, we derive lower bounds for the first eigenvalue of the Laplace--de~Rham operator under natural curvature assumptions, highlighting the distinct roles of holomorphic Ricci and sectional curvatures in Hermitian spectral geometry.

Our approach is based on Bochner type identities adapted to the Strominger--Bismut connection. The balanced condition enables generalizations of differential and integral formulas from the Riemannian to the Hermitian context, and provides new insights into eigenvalue phenomena beyond the K\"ahler realm.

The main theorems can be summarized as follows.

The first one is a Lichnerowicz--Obata type estimate and rigidity statement for the Laplace--de Rham operator on complete balanced Hermitian manifolds.

\begin{theorem}\label{thm1.1}
Let $(M,\omega)$ be a complete balanced Hermitian manifold of complex dimension $n$. If there exists some constant $K>0$ such that the
holomorphic Ricci curvature of the Strominger--Bismut connection satisfies
\begin{equation}\label{1.1x}
\mathcal{R}ic^{SB,\mathbb{C}}(W,\overline{W})\geq(2n-1)K|W|^2\quad\text{for all}\quad W\in\Gamma(M,T^{1,0}M),
\end{equation}
then
\begin{equation}\label{1.2x}
\lambda_1\geq 2nK.
\end{equation}
When equality is achieved in \eqref{1.2x}, one has
\begin{equation}\label{1.2xx}
D = \frac{\pi}{\sqrt{K}}.
\end{equation}
Furthermore, if $\omega$ is K\"ahler, then $(M,\omega)$ is, up to a scaling, isometric to the complex projective line $\mathbb{CP}^1$ endowed with the Fubini-Study metric $\omega_{FS}$.
\end{theorem}

\begin{example}
On $(\mathbb{C}\mathbb{P}^n,\omega_{\mathrm{FS}})$, it is well-known that
\[
Ric(\omega_{\mathrm{FS}})=(n+1)\omega_{\mathrm{FS}}, \qquad
\mathrm{HSC}=2, \qquad\text{and}\qquad
\lambda_1=2(n+1).
\]
\end{example}

This example naturally raises the question of whether the rigidity phenomenon in Theorem~\ref{thm1.1} persists in a more general, possibly non-K\"ahler, setting.

\begin{conjecture}
The rigidity conclusion of Theorem~\ref{thm1.1} that $(M,\omega)$ is isometric to $(\mathbb{CP}^1, \omega_{\mathrm{FS}})$ up to a scaling remains valid without assuming $\omega$ is K\"ahler.
\end{conjecture}

In the compact case, the curvature condition can be further relaxed.

\begin{corollary}\label{cor1.1}
Let $(M,\omega)$ be a compact balanced Hermitian manifold of complex dimension $n$. If there exists some constant $K>0$ such that the holomorphic Ricci curvature of the Strominger--Bismut connection satisfies
\begin{equation}
\mathcal{R}ic^{SB,\mathbb{C}}(U,\overline{U})\geq(2n-1)K|U|^2,
\end{equation}
then
\begin{equation}\label{1.2}
\lambda_1\geq 2nK.
\end{equation}
\end{corollary}

The Li--Yau type first eigenvalue estimate on compact Riemannian manifolds \cite{LY80} can be extended to compact balanced Hermitian manifolds of complex dimension $n\geq3$.

\begin{theorem}\label{thm1.2}
Let $(M,\omega)$ be a compact balanced Hermitian manifold of complex dimension $n\geq3$. If there exists some constant $K\geq0$ such that the holomorphic Ricci curvature of the Strominger--Bismut connection satisfies
\begin{equation}\label{1.12x}
\mathcal{R}ic^{SB,\mathbb{C}}(U,\overline{U})\geq -K|U|^2,
\end{equation}
then there exist positive constants $C_1$ and $C_2$ depending only on $n$, such that
\begin{equation}\label{1.13x}
\lambda_1\geq\frac{C_1}{D^2}\exp(-C_2\sqrt{K}D).
\end{equation}
\end{theorem}

Theorem \ref{thm1.2} admits a further extension to complete balanced Hermitian manifolds under stronger holomorphic Ricci curvature bounds.

\begin{corollary}\label{cor1.2}
Let $(M,\omega)$ be a complete balanced Hermitian manifold of complex dimension $n\geq3$. If there exists some constant $K>0$ such that the holomorphic Ricci curvature satisfies
\begin{equation}\label{1.12}
\mathcal{R}ic^{SB,\mathbb{C}}(W,\overline{W})\geq K|W|^2\quad\text{for all}\quad W\in\Gamma(M,T^{1,0}M),
\end{equation}
then
\begin{equation}\label{1.13}
\lambda_1\geq\frac{4}{(2n+1)e^2D^2}.
\end{equation}
\end{corollary}

Similarly, the Zhong--Yang type first eigenvalue estimate \cite{ZY84} extends to complete balanced Hermitian manifolds.

\begin{theorem}\label{thm1.3}
Let $(M,\omega)$ be a complete balanced Hermitian manifold. If there exists some constant $K>0$ such that the holomorphic Ricci curvature satisfies
\begin{equation}\label{1.14}
\mathcal{R}ic^{SB,\mathbb{C}}(W,\overline{W})\geq K|W|^2\quad\text{for all}\quad W\in\Gamma(M,T^{1,0}M),
\end{equation}
then
\begin{equation}\label{1.15}
\lambda_1\geq \frac{\pi^2}{D^2}.
\end{equation}
\end{theorem}

\begin{corollary}\label{cor1.3}
Let $(M,\omega)$ be a compact balanced Hermitian manifold. If the holomorphic Ricci curvature satisfies
\begin{equation}
\mathcal{R}ic^{SB,\mathbb{C}}(U,\overline{U})\geq 0,
\end{equation}
then
\begin{equation}
\lambda_1\geq \frac{\pi^2}{D^2}.
\end{equation}
\end{corollary}

Finally, considering lower bounds on the holomorphic sectional curvature of the Strominger--Bismut connection, we have:

\begin{theorem}\label{thm1.4}
Let $(M,\omega)$ be a complete balanced Hermitian manifold. If $[\Lambda,\bar{\partial}\omega]\partial u=0$ and
\begin{equation}\label{6.8}
\mathrm{HSC}^{SB}\geq K
\end{equation}
for some constant $K>0$,
then
\begin{equation}\label{6.15}
\lambda_1\geq K.
\end{equation}
\end{theorem}

The first-eigendirection torsion condition $[\Lambda,\bar{\partial}\omega]\partial u=0$ is implied by the geometric assumption that the torsion tensor ${^{SB}}T$ of the Strominger--Bismut connection vanishes in the direction of $U$, namely, ${^{SB}}T_{ij}^kU^i=0$ for all $j$, $k$. Notice that $U+\overline{U}\in T_{\mathbb{R}}M$. Moreover, in the compact case, the holomorphic sectional curvature condition \eqref{6.8} can be further weakened.

\begin{corollary}\label{cor1.5}
Let $(M,\omega)$ be a compact balanced Hermitian manifold. Assume that ${^{SB}}T_{ij}^kU^i=0$ for all $j$, $k$, and at every point where $U\neq0$,
\begin{equation}
\mathrm{HSC}^{SB}(U+\overline{U})\geq K
\end{equation}
for some constant $K>0$. Then
\begin{equation}
\lambda_1\geq K.
\end{equation}
\end{corollary}

The paper is organized as follows.
Section~\ref{sec2} fixes notation and recalls some formulas for the Chern and Strominger--Bismut connections.
Sections~\ref{sec3} to \ref{sec5} establish eigenvalue estimates under holomorphic Ricci curvature bounds.
Section~\ref{sec6} derives bounds in terms of holomorphic sectional curvature.
\section{Preliminaries}
\label{sec2}
\subsection{The Chern connection}
Let $\{z^1,z^2,\cdots,z^n\}$ be the local holomorphic coordinates on a Hermitian manifold $(M,h)$. The Chern connection $^C\nabla$ is the unique Hermitian connection that is compatible with the holomorphic structure. Its coefficients on the holomorphic tangent bundle $T^{1,0}M$ satisfy
\begin{equation}
{^C}\Gamma_{ij}^k=h^{k\bar{l}}\frac{\partial h_{j\bar{l}}}{\partial z^i}\quad\text{and}\quad {^C}\Gamma_{\bar{i}j}^k={^C}\Gamma_{\bar{i}j}^{\bar{k}}={^C}\Gamma_{ij}^{\bar{k}}=0.
\end{equation}

Let $^C\nabla^{1,0}$ and $^C\nabla^{0,1}$ denote $(1,0)$- and $(0,1)$-component of the Chern connection $^C\nabla$, respectively. Then
\[{^C}\nabla={^C}\nabla^{1,0}+{^C}\nabla^{0,1}\quad\text{with}\quad{^C}\nabla^{0,1}=\bar{\partial}.\]

The torsion tensor ${^C}T$ of the Chern connection ${^C}\nabla$ is
\begin{equation}
{^C}T_{ij}^k={^C}\Gamma_{ij}^k-{^C}\Gamma_{ji}^k=h^{k\bar{l}}\big(\frac{\partial h_{j\bar{l}}}{\partial z^i}-\frac{\partial h_{i\bar{l}}}{\partial z^j}\big).
\end{equation}

The curvature tensor $\Theta$ of the Chern connection ${^C}\nabla$ on the Hermitian holomorphic tangent bundle $(T^{1,0}M,h)$ is given by
\begin{equation}
\Theta_{i\bar{j}k\bar{l}}=h_{p\bar{l}}\Theta_{i\bar{j}k}^p=-h_{p\bar{l}}\frac{\partial{^C}\Gamma_{ik}^p}{\partial\bar{z}^j}=-\frac{\partial^2h_{k\bar{l}}}{\partial z^i\partial\bar{z}^j}+h^{p\bar{q}}\frac{\partial h_{p\bar{l}}}{\partial\bar{z}^j}\frac{\partial h_{k\bar{q}}}{\partial z^i},
\end{equation}
and the first Chern--Ricci curvature by
\begin{equation}
\Theta^{(1)}=\sqrt{-1}\, \Theta^{(1)}_{i\bar{j}}dz^i\wedge d\bar{z}^j,
\end{equation}
where
\[\Theta^{(1)}_{i\bar{j}}=h^{k\bar{l}}\Theta_{i\bar{j}k\bar{l}}=-\frac{\partial^2\log\det(h_{k\bar{l}})}{\partial z^i\partial\bar{z}^j}.\]

\subsection{The Strominger--Bismut connection.}
Let $g$ denote the background Riemannian metric and $J$ the complex structure, satisfying
\[
g(X,Y) = g(JX,JY), \quad \omega(X,Y) = g(JX,Y)
\]
for any $X,Y \in \Gamma(M,T_{\mathbb{R}}M)$.

Let \(g^{\mathbb C}\) denote the \(\mathbb C\)-bilinear extension of
\(g\) to $$T_{\mathbb{C}}M=T_{\mathbb{R}}M\otimes\mathbb{C}=T^{1,0}M\oplus T^{0,1}M.$$
It is known that
\[
g^{\mathbb C}(\frac{\partial}{\partial z^i},
\frac{\partial}{\partial\bar z^j})=h(\frac{\partial}{\partial z^i},
\frac{\partial}{\partial z^j})=h_{i\bar j},
\quad\text{and}\quad
g^{\mathbb C}(\frac{\partial}{\partial z^i},
\frac{\partial}{\partial z^j})=0.
\]

The Strominger--Bismut connection $^{SB}\nabla$ is the unique Hermitian connection whose torsion ${^{SB}}T$ is totally skew-symmetric, namely, $${^{SB}}\nabla g=0,\quad{^{SB}}\nabla J=0\quad\text{and}\quad{^{SB}}T\in\Gamma(M,\wedge^3T^{*}_{\mathbb{R}}M),$$
where ${^{SB}}T$, regarded as a 3--form, is given by
\[^{SB}T(X,Y,Z):=g({^{SB}}\nabla_XY-{^{SB}}\nabla_YX-[X,Y],Z),\quad X,Y,Z\in \Gamma(M,T_{\mathbb{R}}M).\]

The coefficients of the Strominger--Bismut connection ${^{SB}}\nabla$  on the holomorphic tangent bundle $T^{1,0}M$ satisfy
\begin{equation}
{^{SB}}\Gamma_{ij}^k=h^{k\bar{l}}\frac{\partial h_{i\bar{l}}}{\partial z^j}={^{C}}\Gamma_{ji}^k,
\end{equation}
\begin{equation}
{^{SB}}\Gamma_{\bar{i}j}^k=h^{k\bar{l}}\big(\frac{\partial h_{j\bar{l}}}{\partial \bar{z}^i}-\frac{\partial h_{j\bar{i}}}{\partial \bar{z}^l}\big),
\end{equation}
and
\begin{equation}
^{SB}\Gamma_{i\bar{j}}^k={^{SB}}\Gamma_{ij}^{\bar{k}}=0,
\end{equation}

We can decompose the Strominger--Bismut connection ${^{SB}}\nabla$ into its $(1,0)$-component ${^{SB}}\nabla^{1,0}$ and $(0,1)$-component ${^{SB}}\nabla^{0,1}$, i.e.,
\begin{equation}
{^{SB}}\nabla={^{SB}}\nabla^{1,0}+{^{SB}}\nabla^{0,1}.
\end{equation}

Given any function $f\in C^\infty(M,\mathbb{R})$, it is clear that
\begin{equation}\label{3.20}
{^{SB}}\nabla^{1,0}\partial f=\big(\frac{\partial^2f}{\partial z^i\partial z^j}-{^{SB}}\Gamma_{ij}^k\frac{\partial f}{\partial z^k}\big)dz^i\otimes dz^j,
\end{equation}
and
\begin{equation}\label{3.21}
{^{SB}}\nabla^{0,1}\partial f=\big(\frac{\partial^2f}{\partial \bar{z}^j\partial z^i}-{^{SB}}\Gamma_{\bar{j}i}^k\frac{\partial f}{\partial z^k}\big)d\bar{z}^j\otimes dz^i.
\end{equation}

For simplicity, we denote
\begin{equation}\label{6.4}
t_{ij}:=\frac{\partial^2u}{\partial z^i\partial z^j}-{^{SB}}\Gamma_{ij}^k\frac{\partial u}{\partial z^k},\qquad t_{\bar{i}\bar{j}}=\overline{t_{ij}},
\end{equation}
and
\begin{equation}\label{6.5}
s_{\bar{j}i}=\frac{\partial^2u}{\partial \bar{z}^j\partial z^i}-{^{SB}}\Gamma_{\bar{j}i}^k\frac{\partial u}{\partial z^k},\qquad s_{i\bar{j}}=\overline{s_{\bar{i}j}}
\end{equation}
throughout this paper.

The torsion tensor ${^{SB}}T$ of the Strominger--Bismut connection ${^{SB}}\nabla$ is
\begin{equation}
{^{SB}}T_{ij}^k={^{SB}}\Gamma_{ij}^k-{^{SB}}\Gamma_{ji}^k=h^{k\bar{l}}\big(\frac{\partial h_{i\bar{l}}}{\partial z^j}-\frac{\partial h_{j\bar{l}}}{\partial z^i}\big)={^{C}}T_{ji}^k.
\end{equation}

The real curvature tensor of the Strominger--Bismut connection $^{SB}\nabla$ on the underlying Riemannian $2n$--manifold $(M,g,J)$ is defined as
\[R^{SB,\mathbb{R}}(X,Y,Z,W)=g({^{SB}}\nabla_X{^{SB}}\nabla_YZ-{^{SB}}\nabla_Y{^{SB}}\nabla_XZ-{^{SB}}\nabla_{[X,Y]}Z,W),\]
where $X$, $Y$, $Z$, $W\in \Gamma(M,T_{\mathbb{R}}M)$.

The holomorphic sectional curvature $\mathrm{HSC}^{SB}$ of the Strominger--Bismut connection $^{SB}\nabla$ in the direction of $X\in T_{\mathbb{R}}M$ is defined by
\begin{equation}\label{1.16}
\mathrm{HSC}^{SB}(X)=\frac{R^{SB,\mathbb{R}}(JX,X,X,JX)}{g(X,X)^2}.
\end{equation}

The real Ricci curvature tensor of the Strominger--Bismut connection $^{SB}\nabla$ on $(M,g,J)$ is defined by
\[\mathcal{R}ic^{SB,\mathbb{R}}(X,Y)=g^{il}R^{SB,\mathbb{R}}\big(\frac{\partial}{\partial x^i},X,Y,\frac{\partial}{\partial x^l}\big)\]
for $X,Y\in \Gamma(M,T_{\mathbb{R}}M)$.

We denote by $R^{SB,\mathbb{C}}(X,Y,Z,W)$ and $\mathcal{R}ic^{SB,\mathbb{C}}(Z,W)$, where $X$, $Y$, $Z$, $W\in \Gamma(M,T_{\mathbb{C}}M)$, the $\mathbb{C}$-linear complexified curvature tensor and Ricci tensor of the Strominger--Bismut connection $^{SB}\nabla$, respectively. Following the definition in \cite{Yang25a+}, the holomorphic Ricci curvature of the Strominger--Bismut connection $^{SB}\nabla$ is given by
\begin{equation}\label{1.1}
\mathcal{R}ic^{SB,\mathbb{C}}(W,\overline{W})=h^{i\bar{l}}R^{SB,\mathbb{C}}\big(\frac{\partial}{\partial z^i},\overline{W},W,\frac{\partial}{\partial \bar{z}^l}\big)
\end{equation}
for $W\in\Gamma(M,T^{1,0}M)$, serving as the complex analogue of the real Ricci curvature.

The first, second, third and fourth Strominger--Bismut--Ricci curvatures of the Hermitian holomorphic tangent bundle $(T^{1,0}M,h)$ are defined as
\begin{equation}
Ric^{SB(1)}=\sqrt{-1}\, R^{SB(1)}_{i\bar{j}}dz^i\wedge d\bar{z}^{j}\quad\text{with}\quad R^{SB(1)}_{i\bar{j}}=h^{k\bar{l}}R^{SB,\mathbb{C}}_{i\bar{j}k\bar{l}},
\end{equation}
\begin{equation}
Ric^{SB(2)}=\sqrt{-1}\, R^{SB(2)}_{i\bar{j}}dz^i\wedge d\bar{z}^{j}\quad\text{with}\quad R^{SB(2)}_{i\bar{j}}=h^{k\bar{l}}R^{SB,\mathbb{C}}_{k\bar{l}i\bar{j}},
\end{equation}
\begin{equation}
Ric^{SB(3)}=\sqrt{-1}\, R^{SB(3)}_{i\bar{j}}dz^i\wedge d\bar{z}^{j}\quad\text{with}\quad R^{SB(3)}_{i\bar{j}}=h^{k\bar{l}}R^{SB,\mathbb{C}}_{i\bar{l}k\bar{j}},
\end{equation}
and
\begin{equation}\label{2.9}
Ric^{SB(4)}=\sqrt{-1}\, R^{SB(4)}_{i\bar{j}}dz^i\wedge d\bar{z}^{j}\quad\text{with}\quad R^{SB(4)}_{i\bar{j}}=h^{k\bar{l}}R^{SB,\mathbb{C}}_{k\bar{j}i\bar{l}},
\end{equation}
respectively, where $R^{SB,\mathbb{C}}_{i\bar{j}k\bar{l}}=R^{SB,\mathbb{C}}(\frac{\partial}{\partial z^i},\frac{\partial}{\partial \bar{z}^j},\frac{\partial}{\partial z^k},\frac{\partial}{\partial \bar{z}^l})$.

It follows from \eqref{1.1} and \eqref{2.9} that
\begin{equation}\label{2.10}
\mathcal{R}ic^{SB,\mathbb{C}}(W,\overline{W})=R_{i\bar{j}}^{SB(4)}W^i\overline{W^j}=\langle Ric^{SB(4)},\sqrt{-1}W^\flat\wedge\overline{W}^\flat),
\end{equation}
where $W=W^i\frac{\partial}{\partial z^i}\in\Gamma(M,T^{1,0}M)$ and $W^\flat=h_{i\bar{j}}\overline{W^j}dz^i$ is its dual $(1,0)$-form.

\subsection{Some operators and identities}
For any tensors (or forms) $\alpha$ and $\beta$ of the same bidegree, we denote by $\langle \alpha,\beta \rangle$ their pointwise inner product with $|\alpha|^2=\langle\alpha,\alpha\rangle$.

Define the operator $(\cdot,\cdot)$ by
\begin{equation}
(\alpha,\beta):=\int_M\langle\alpha,\beta\rangle\frac{\omega^n}{n!}\quad\text{with}\quad\|\alpha\|^2:=(\alpha,\alpha).
\end{equation}

Let $\{\cdot,\cdot\}$ be the pointwise inner product on $\wedge^rT^*_{\mathbb{C}}M$-valued complex 1-forms. For example,
\[\{{^{SB}}\nabla^{1,0}\partial u,\partial u\}=h^{j\bar{k}}t_{ij}\frac{\partial u}{\partial \bar{z}^k}dz^i\quad\text{and}\quad\{\partial u,{^{SB}}\nabla^{0,1}\partial u\}=h^{j\bar{k}}\frac{\partial u}{\partial z^j}s_{i\bar{k}}dz^i.\]

The curvature relation between the Strominger--Bismut connection and the Chern connection on a balanced Hermitian manifold can be directly deduced from \cite[Corollary 1.8]{WY25}.
\begin{lemma}[\cite{WY25}]\label{lem6.1}
Let $(M,\omega)$ be a complete balanced Hermitian manifold. The complexification of the curvature tensor of the Strominger--Bismut connection is given by
\begin{equation}\label{6.0}
R^{SB,\mathbb{C}}_{i\bar{j}k\bar{l}}=\Theta_{i\bar{l}k\bar{j}}+\Theta_{k\bar{j}i\bar{l}}-\Theta_{i\bar{j}k\bar{l}}+h_{p\bar{q}}{^{C}}T_{ik}^p\overline{{^{C}}T_{jl}^q}-h^{p\bar{q}}h_{m\bar{l}}h_{k\bar{s}}{^{C}}T_{ip}^m\overline{{^{C}}T_{jq}^s}
\end{equation}
and the Ricci curvatures are
\begin{equation}\label{6.1}
Ric^{SB(1)}=\Theta^{(1)}=\Theta^{(3)}=\Theta^{(4)},
\end{equation}
\begin{equation}\label{6.2}
Ric^{SB(2)}=\Theta^{(1)}+\sqrt{-1}\, \Lambda\partial\bar{\partial}\omega-\sqrt{-1}\, T\circ\overline{T},
\end{equation}
where
\begin{equation}\label{6.7}
T\circ\overline{T}=h^{p\bar{q}}h^{s\bar{t}}h_{k\bar{j}}h_{i\bar{l}}{^C}T_{sp}^k\overline{{^C}T_{tq}^l}dz^i\wedge d\bar{z}^j.
\end{equation}
and
\begin{equation}\label{6.3}
Ric^{SB(3)}=Ric^{SB(4)}=\Theta^{(1)}-\sqrt{-1}\, \Lambda\partial\bar{\partial}\omega.
\end{equation}
\end{lemma}

We obtain the following corollary concerning the holomorphic Ricci curvature of the Strominger--Bismut connection.
\begin{corollary}\label{cor6.2}
Let $(M,\omega)$ be a complete Hermitian manifold. The holomorphic Ricci curvature of the Strominger--Bismut connection is non-negative, i.e.,$$\mathcal{R}ic^{SB,\mathbb{C}}(W,\overline{W}) \geq 0\quad \text{for all}\quad W \in \Gamma(M,T^{1,0}M)$$ if and only if $$Ric^{SB(3)}+Ric^{SB(4)} \geq 0.$$

In particular, if $\omega$ is balanced, one has
\begin{equation}\label{6.6}
\mathcal{R}ic^{SB,\mathbb{C}}(W,\overline{W}) = -(2\sqrt{-1}Ric^{SB(1)}-\sqrt{-1}Ric^{SB(2)}+\, T \circ \overline{T})(W,\overline{W})
\end{equation}
for $W \in \Gamma(M,T^{1,0}M)$.
\end{corollary}

\section{Lichnerowicz--Obata type estimates}
\label{sec3}
In this section, we establish Theorem~\ref{thm1.1} and Corollary~\ref{cor1.1}, which yield Lichnerowicz--Obata type lower bounds for the first eigenvalue of the Laplace--de~Rham operator on complete balanced Hermitian manifolds.
\begin{lemma}\label{lem3.1}
Given $f\in C^\infty(M,\mathbb{R})$ on a compact balanced Hermitian manifold $(M,\omega)$, one has
\begin{equation}\label{3.0}
\frac{1}{2}\Delta_df=\Delta_{\partial}f=\Delta_{\bar{\partial}}f=-\mathrm{tr}_\omega(\sqrt{-1}\, \partial\bar{\partial}f).
\end{equation}
\end{lemma}

\begin{proof}
It is well-known that
\begin{equation}\label{3.3}
\Delta_df=d^*df=(\partial^*+\bar{\partial}^*)(\partial+\bar{\partial})f=\partial^*\partial f+\bar{\partial}^*\bar{\partial}f=\Delta_{\partial}f+\Delta_{\bar{\partial}}f.
\end{equation}

For any $F\in C^\infty(M,\mathbb{R})$, we have
\begin{eqnarray}\label{3.5}
(\Delta_\partial f,F)&=&(\partial^*\partial f,F)\notag\\
&=&(\partial f,\partial F)\notag\\
&=&\int_M\mathrm{tr}_\omega(\sqrt{-1}\, \partial f\wedge\bar{\partial} F)\frac{\omega^n}{n!}\notag\\
&=&\int_M\sqrt{-1}\, \partial f\wedge\bar{\partial} F\wedge\frac{\omega^{n-1}}{(n-1)!}\notag\\
&=&\int_M\bar{\partial}(-\sqrt{-1}\, F\partial f\wedge\frac{\omega^{n-1}}{(n-1)!})\notag\\
&&+\int_M\sqrt{-1}\, F\bar{\partial}\partial f\wedge\frac{\omega^{n-1}}{(n-1)!}\notag\\
&&-\int_M\sqrt{-1}\, F\partial f\wedge\frac{\bar{\partial}\omega^{n-1}}{(n-1)!}\notag\\
&=&-\int_MF\mathrm{tr}_\omega(\sqrt{-1}\, \partial\bar{\partial}f)\frac{\omega^n}{n!},
\end{eqnarray}
where we used the Stokes' formula and the condition of $\bar{\partial}\omega^{n-1}=0$.
Moreover, we obtain that
\begin{eqnarray}\label{3.6}
(\Delta_{\bar{\partial}}f,F)&=&\overline{(\Delta_{\partial}f,F)}=-\int_MF\mathrm{tr}_\omega(\sqrt{-1}\, \partial\bar{\partial}f)\frac{\omega^n}{n!}.
\end{eqnarray}

\eqref{3.0} follows by \eqref{3.3}, \eqref{3.5} and \eqref{3.6}.
\end{proof}

\begin{lemma}\label{lem3.2}
Given a $(0,1)$-form $\varphi$ on a compact balanced Hermitian manifold $(M,\omega)$, one has
\begin{equation}\label{3.7}
\bar{\partial}^*\varphi=-\sqrt{-1}\, \Lambda\partial\varphi,
\end{equation}
and
\begin{equation}\label{3.7x}
\partial^*\bar{\varphi}=\sqrt{-1}\, \Lambda\bar{\partial}\bar{\varphi}.
\end{equation}
\end{lemma}

\begin{proof}
First of all, we claim that
\begin{equation}\label{3.8}
\sqrt{-1}\, \langle\varphi,\partial^*\omega\rangle=\sqrt{-1}\, \Lambda\partial\varphi+\bar{\partial}^*\varphi
\end{equation}
on compact Hermitian manifolds.

For any $F\in C^\infty(M,\mathbb{R})$, we have
\begin{eqnarray*}
(\sqrt{-1}\, \langle\varphi,\partial^*\omega\rangle,F)&=&(\sqrt{-1}\, F\varphi,\partial^*\omega)\notag\\
&=&(\partial(\sqrt{-1}\, F\varphi),\omega)\notag\\
&=&(\sqrt{-1}\, \partial F\wedge\varphi,\omega)+(F\sqrt{-1}\, \partial\varphi,\omega)\notag\\
&=&\int_M\mathrm{tr}_\omega(\sqrt{-1}\, \partial F\wedge\varphi)\frac{\omega^n}{n!}+\int_MF\mathrm{tr}_\omega(\sqrt{-1}\, \partial\varphi)\frac{\omega^n}{n!}\notag\\
&=&(\partial F,\varphi)+(\mathrm{tr}_\omega(\sqrt{-1}\, \partial\varphi),F)\notag\\
&=&(\bar{\partial}^*\varphi+\sqrt{-1}\, \Lambda\partial\varphi,F).
\end{eqnarray*}
This proves \eqref{3.8}.

Since $\partial^*\omega=0$, we get \eqref{3.7} from \eqref{3.8}. \eqref{3.7x} follows by
\[\partial^*\bar{\varphi}=\overline{\bar{\partial}^*\varphi}=\overline{-\sqrt{-1}\, \Lambda\partial\varphi}=\sqrt{-1}\, \Lambda\bar{\partial}\bar{\varphi}.\]
\end{proof}

\begin{proposition}\label{prp3.2}
On a compact balanced Hermitian manifold $(M,\omega)$, one has
\begin{equation}\label{3.28}
h^{i\bar{j}}s_{i\bar{j}}=\mathrm{tr}_\omega(\sqrt{-1}\, \partial\bar{\partial}u)=-\frac{\lambda_1}{2}u,
\end{equation}
where $s_{i\bar{j}}$ is defined in \eqref{6.5}.
\end{proposition}

\begin{proof}
It is known that (see, e.g., \cite[Lemma 3.3]{LY17}, \cite[(3.5)]{Yang25b+})
\begin{eqnarray}\label{3.29}
\bar{\partial}^*\omega&=&\sqrt{-1}\, \Lambda(\partial\omega)\notag\\
&=&h^{i\bar{j}}\iota_{\bar{j}}\iota_{i}\big(\sqrt{-1}\, \frac{\partial h_{p\bar{q}}}{\partial z^k}dz^k\wedge dz^p\wedge d\bar{z}^q\big)\notag\\
&=&\sqrt{-1}\, h^{i\bar{j}}\big(\frac{\partial h_{i\bar{j}}}{\partial z^k}-\frac{\partial h_{k\bar{j}}}{\partial z^i}\big)dz^k.
\end{eqnarray}

Since $\bar{\partial}^*\omega=0$, we obtain that
\[h^{i\bar{j}}\big(\frac{\partial h_{i\bar{j}}}{\partial z^k}-\frac{\partial h_{k\bar{j}}}{\partial z^i}\big)=0.\]

Moreover, we have
\begin{eqnarray}
h^{i\bar{j}}s_{i\bar{j}}&=&\mathrm{tr}_\omega(\sqrt{-1}\, \partial\bar{\partial}u)-h^{i\bar{j}}\overline{{^{SB}}\Gamma_{\bar{i}j}^k}\frac{\partial u}{\partial\bar{z}^k}\notag\\
&=&\mathrm{tr}_\omega(\sqrt{-1}\, \partial\bar{\partial}u)-h^{i\bar{j}}h^{l\bar{k}}\big(\frac{\partial h_{l\bar{j}}}{\partial z^i}-\frac{\partial h_{i\bar{j}}}{\partial z^l}\big)\frac{\partial u}{\partial\bar{z}^k}\notag\\
&=&\mathrm{tr}_\omega(\sqrt{-1}\, \partial\bar{\partial}u).
\end{eqnarray}
Combining with Lemma \ref{lem3.1} and \eqref{1}, we get \eqref{3.28}.
\end{proof}

\begin{proposition}\label{prp3.1}
On a compact balanced Hermitian manifold $(M,\omega)$, one has the Bochner type formula that
\begin{equation}\label{3.12}
\Delta_{\bar{\partial}}|\partial u|^2=-\mathcal{R}ic^{SB,\mathbb{C}}(U,\overline{U})+\lambda_1|\partial u|^2-|{^{SB}}\nabla^{1,0}\partial u|^2-|{^{SB}}\nabla^{0,1}\partial u|^2.
\end{equation}
\end{proposition}

\begin{proof}
Note that
\begin{eqnarray}\label{3.13}
\Delta_{\bar{\partial}}|\partial u|^2&=&\bar{\partial}^*\bar{\partial}|\partial u|^2\notag\\
&=&\bar{\partial}^*(\{{^{SB}}\nabla^{0,1}\partial u,\partial u\}+\{\partial u,{^{SB}}\nabla^{1,0}\partial u\})\notag\\
&=&-\sqrt{-1}\, \Lambda(\partial\{{^{SB}}\nabla^{0,1}\partial u,\partial u\}+\partial\{\partial u,{^{SB}}\nabla^{1,0}\partial u\})\notag\\
&=&-\sqrt{-1}\, \Lambda(\partial\{{^{SB}}\nabla^{0,1}\partial u,\partial u\}+\bar{\partial}\{{^{SB}}\nabla^{1,0}\partial u,\partial u\})\notag\\
&&+\sqrt{-1}\, \Lambda(\bar{\partial}\{{^{SB}}\nabla^{1,0}\partial u,\partial u\}-\partial\{\partial u,{^{SB}}\nabla^{1,0}\partial u\}),
\end{eqnarray}
where we used \eqref{3.7} in the third equality.

In the following, we deal with the right-hand-side of \eqref{3.13}.

Using \eqref{3.20} and \eqref{3.21}, we obtain that
\begin{eqnarray}\label{3.14}
&&-\sqrt{-1}\, \Lambda(\partial\{{^{SB}}\nabla^{0,1}\partial u,\partial u\}+\bar{\partial}\{{^{SB}}\nabla^{1,0}\partial u,\partial u\})\notag\\
&=&\sqrt{-1}\, \Lambda\{h^{b\bar{l}}R_{i\bar{j}k\bar{l}}^{SB,\mathbb{C}}\frac{\partial u}{\partial z^b}dz^i\wedge d\bar{z}^j\otimes dz^k,\frac{\partial u}{\partial z^a}dz^a\}\notag\\
&&+\sqrt{-1}\, \Lambda\{{^{SB}}\nabla^{1,0}\partial u,{^{SB}}\nabla^{1,0}\partial u\}\notag\\
&&+\sqrt{-1}\, \Lambda\{{^{SB}}\nabla^{0,1}\partial u,{^{SB}}\nabla^{0,1}\partial u\}\notag\\
&=&h^{p\bar{q}}\iota_{\bar{q}}\iota_p(h^{k\bar{a}}h^{b\bar{l}}R^{SB,\mathbb{C}}_{i\bar{j}k\bar{l}}\frac{\partial u}{\partial \bar{z}^a}\frac{\partial u}{\partial z^b}dz^i\wedge d\bar{z}^j)\notag\\
&&+h^{p\bar{q}}\iota_{\bar{q}}\iota_p(h^{k\bar{l}}t_{ik}t_{\bar{j}\bar{l}}dz^i\wedge d\bar{z}^j)\notag\\
&&+h^{p\bar{q}}\iota_{\bar{q}}\iota_p(h^{k\bar{l}}s_{\bar{j}k}s_{i\bar{l}}d\bar{z}^j\wedge dz^i)\notag\\
&=&h^{k\bar{a}}h^{b\bar{l}}R_{k\bar{l}}^{SB(2)}\frac{\partial u}{\partial \bar{z}^a}\frac{\partial u}{\partial z^b}+h^{i\bar{j}}h^{k\bar{l}}t_{ik}t_{\bar{j}\bar{l}}-h^{i\bar{j}}h^{k\bar{l}}s_{\bar{j}k}s_{i\bar{l}}\notag\\
&=&-\sqrt{-1}\, Ric^{SB(2)}(U,\overline{U})+|{^{SB}}\nabla^{1,0}\partial u|^2-|{^{SB}}\nabla^{0,1}\partial u|^2.
\end{eqnarray}

By \eqref{3.20}, we have
\begin{eqnarray}\label{3.15}
&&\sqrt{-1}\, \Lambda\bar{\partial}\{{^{SB}}\nabla^{1,0}\partial u,\partial u\}\notag\\
&=&h^{p\bar{q}}\iota_{\bar{q}}\iota_p\big(\frac{\partial}{\partial\bar{z}^l}\big(h^{j\bar{k}}t_{ij}\frac{\partial u}{\partial \bar{z}^k}\big)d\bar{z}^l\wedge dz^i\big)\notag\\
&=&-h^{p\bar{q}}\frac{\partial}{\partial\bar{z}^q}\big(h^{i\bar{k}}t_{pi}\frac{\partial u}{\partial \bar{z}^k}\big)\notag\\
&=&h^{p\bar{q}}h^{i\bar{m}}h^{r\bar{k}}\frac{\partial h_{r\bar{m}}}{\partial\bar{z}^q}t_{pi}\frac{\partial u}{\partial \bar{z}^k}-h^{p\bar{q}}h^{i\bar{k}}t_{pi}\frac{\partial^2 u}{\partial \bar{z}^k\partial\bar{z}^q}\notag\\
&&-h^{p\bar{q}}h^{i\bar{k}}(\frac{\partial^3u}{\partial z^i\partial z^p\partial\bar{z}^q}-\frac{\partial{^{SB}}\Gamma_{pi}^s}{\partial\bar{z}^q}\frac{\partial u}{\partial z^s}-{^{SB}}\Gamma_{pi}^s\frac{\partial^2 u}{\partial z^s\partial\bar{z}^q})\frac{\partial u}{\partial \bar{z}^k}.
\end{eqnarray}

Since ${^{SB}}T_{ij}^k={^{C}}T_{ji}^k$ and is anti--symmetric in $i$ and $j$, it follows that
\begin{eqnarray}\label{3.16}
&&h^{p\bar{q}}h^{i\bar{m}}h^{r\bar{k}}\frac{\partial h_{r\bar{m}}}{\partial\bar{z}^q}t_{pi}\frac{\partial u}{\partial \bar{z}^k}-h^{p\bar{q}}h^{i\bar{k}}t_{pi}\frac{\partial^2 u}{\partial \bar{z}^k\partial\bar{z}^q}\notag\\
&=&-h^{p\bar{q}}h^{i\bar{m}}t_{pi}\big(\frac{\partial^2u}{\partial \bar{z}^m\partial \bar{z}^q}-\overline{{^{SB}}\Gamma_{qm}^k}\frac{\partial u}{\partial z^k}\big)\notag\\
&&-h^{p\bar{q}}h^{i\bar{m}}t_{pi}\overline{{^{SB}}T_{qm}^k}\frac{\partial u}{\partial z^k}\notag\\
&=&-h^{p\bar{q}}h^{i\bar{m}}t_{pi}\overline{t_{qm}}+\frac{1}{2}h^{p\bar{q}}h^{i\bar{m}}{^{C}}T_{ip}^l\frac{\partial u}{\partial z^l}\overline{{^{C}}T_{mq}^k}\frac{\partial u}{\partial z^k}\notag\\
&=&-|{^{SB}}\nabla^{1,0}\partial u|^2+\frac{1}{2}\, T\circ\overline{T}(U,\overline{U}).
\end{eqnarray}

Moreover, we get
\begin{eqnarray}\label{3.17}
&&-h^{p\bar{q}}h^{i\bar{k}}(\frac{\partial^3u}{\partial z^i\partial z^p\partial\bar{z}^q}-\frac{\partial{^{SB}}\Gamma_{pi}^s}{\partial\bar{z}^q}\frac{\partial u}{\partial z^s}-{^{SB}}\Gamma_{pi}^s\frac{\partial^2 u}{\partial z^s\partial\bar{z}^q})\frac{\partial u}{\partial \bar{z}^k}\notag\\
&=&-h^{p\bar{q}}h^{i\bar{k}}(\frac{\partial^3u}{\partial z^i\partial z^p\partial\bar{z}^q}-\frac{\partial{^{C}}\Gamma_{ip}^s}{\partial\bar{z}^q}\frac{\partial u}{\partial z^s}-{^{SB}}\Gamma_{pi}^s\frac{\partial^2 u}{\partial z^s\partial\bar{z}^q})\frac{\partial u}{\partial \bar{z}^k}\notag\\
&=&-h^{i\bar{k}}\frac{\partial}{\partial z^i}(h^{p\bar{q}}\frac{\partial^2u}{\partial z^p\partial\bar{z}^q})\frac{\partial u}{\partial \bar{z}^k}+h^{i\bar{k}}\frac{\partial h^{p\bar{q}}}{\partial z^i}\frac{\partial^2u}{\partial z^p\partial\bar{z}^q}\frac{\partial u}{\partial \bar{z}^k}\notag\\
&&-h^{p\bar{q}}h^{i\bar{k}}\Theta_{i\bar{q}p}^s\frac{\partial u}{\partial z^s}\frac{\partial u}{\partial \bar{z}^k}+h^{p\bar{q}}h^{i\bar{k}}{^{SB}}\Gamma_{pi}^s\frac{\partial^2 u}{\partial z^s\partial\bar{z}^q}\frac{\partial u}{\partial \bar{z}^k}\notag\\
&=&h^{i\bar{k}}\frac{\partial}{\partial z^i}(\frac{1}{2}\Delta_du)\frac{\partial u}{\partial \bar{z}^k}-h^{i\bar{k}}h^{p\bar{q}}{^{SB}}\Gamma_{pi}^s\frac{\partial^2u}{\partial z^s\partial\bar{z}^q}\frac{\partial u}{\partial \bar{z}^k}\notag\\
&&-h^{i\bar{k}}h^{s\bar{q}}\Theta_{i\bar{q}}^{(3)}\frac{\partial u}{\partial z^s}\frac{\partial u}{\partial \bar{z}^k}+h^{p\bar{q}}h^{i\bar{k}}{^{SB}}\Gamma_{pi}^s\frac{\partial^2 u}{\partial z^s\partial\bar{z}^q}\frac{\partial u}{\partial \bar{z}^k}\notag\\
&=&\frac{\lambda_1}{2}|\partial u|^2+\sqrt{-1}\, Ric^{SB(1)}(U,\overline{U}),
\end{eqnarray}
where we used Lemma \ref{lem3.1} in the third equality and \eqref{6.1} in the last.

Applying \eqref{3.16} and \eqref{3.17} to \eqref{3.15}, we have
\begin{eqnarray}\label{3.18}
\sqrt{-1}\, \Lambda\bar{\partial}\{{^{SB}}\nabla^{1,0}\partial u,\partial u\}&=&\sqrt{-1}\, Ric^{SB(1)}(U,\overline{U})+\frac{1}{2}\, T\circ\overline{T}(U,\overline{U})\notag\\
&&+\frac{\lambda_1}{2}|\partial u|^2-|{^{SB}}\nabla^{1,0}\partial u|^2.
\end{eqnarray}

By \eqref{3.18}, we get
\begin{eqnarray}\label{3.19}
-\sqrt{-1}\, \Lambda\partial\{\partial u,{^{SB}}\nabla^{1,0}\partial u\})&=&\overline{\sqrt{-1}\, \Lambda\bar{\partial}\{{^{SB}}\nabla^{1,0}\partial u,\partial u\}}\notag\\
&=&-\sqrt{-1}\, \overline{Ric^{SB(1)}}(\overline{U},U)+\frac{1}{2}\, \overline{T\circ\overline{T}}(\overline{U},U)\notag\\
&&+\frac{\lambda_1}{2}|\bar{\partial} u|^2-|\overline{{^{SB}}\nabla^{1,0}}\bar{\partial} u|^2\notag\\
&=&\sqrt{-1}\, Ric^{SB(1)}(U,\overline{U})+\frac{1}{2}\, T\circ\overline{T}(U,\overline{U})\notag\\
&&+\frac{\lambda_1}{2}|\partial u|^2-|{^{SB}}\nabla^{1,0}\partial u|^2.
\end{eqnarray}

Applying \eqref{3.14}, \eqref{3.18} and \eqref{3.19} to \eqref{3.13} and using \eqref{6.6}, we obtain \eqref{3.12}.
\end{proof}

In the remainder of this section, it will suffice to establish Theorem \ref{thm1.1}, since the proof of Corollary \ref{cor1.1} is entirely contained therein.
\\\\\textbf{Proof of Theorem \ref{thm1.1}.} As $K>0$, \cite[Theorem 1.5]{Yang25a+} ensures that \eqref{1.1x} implies that
$M$ is compact and
\begin{equation}\label{3.30}
D\leq\frac{\pi}{\sqrt{K}}
\end{equation}

Define $$Q=|\partial u|^2+\frac{\lambda_1}{4n}u^2.$$

By Lemma \ref{lem3.1} and \eqref{3.12}, we get
\begin{eqnarray}\label{3.24}
\Delta_{\bar{\partial}}Q&=&\Delta_{\bar{\partial}}|\partial u|^2-\frac{\lambda_1}{4n}\mathrm{tr}_\omega(\sqrt{-1}\, \partial\bar{\partial}u^2)\notag\\
&=&-\mathcal{R}ic^{SB,\mathbb{C}}(U,\overline{U})+\frac{2n-1}{2n}\lambda_1|\partial u|^2-|{^{SB}}\nabla^{1,0}\partial u|^2\notag\\
&&-|{^{SB}}\nabla^{0,1}\partial u|^2-\frac{\lambda_1u}{2n}\mathrm{tr}_\omega(\sqrt{-1}\, \partial\bar{\partial}u).
\end{eqnarray}

Cauchy-Schwartz inequality and \eqref{3.28} show that
\begin{equation}\label{3.25}
-|{^{SB}}\nabla^{0,1}\partial u|^2\leq-\frac{(h^{i\bar{j}}s_{i\bar{j}})^2}{n}=\frac{\lambda_1u}{2n}\mathrm{tr}_\omega(\sqrt{-1}\, \partial\bar{\partial}u).
\end{equation}

Applying \eqref{1.1x} and \eqref{3.25} to \eqref{3.24}, we obtain that
\begin{equation}\label{3.27}
\Delta_{\bar{\partial}}Q\leq\frac{2n-1}{2n}(\lambda_1-2nK)|\partial u|^2.
\end{equation}
Integrating \eqref{3.27} gives
\[0=\int_M\Delta_{\bar{\partial}}Q\frac{\omega^n}{n!}\leq\frac{2n-1}{2n}(\lambda_1-2nK)\int_M|\partial u|^2\frac{\omega^n}{n!},\]
We conclude \eqref{1.2x} that $\lambda_1\geq 2nK$.

In the following, we consider the case of $\lambda_1=2nK$.

By Lemma \ref{lem3.1} and \eqref{3.27}, we have
\[\mathrm{tr}_\omega(\sqrt{-1}\, \partial\bar{\partial}Q)=-\Delta_{\bar{\partial}}Q\geq0.\]
The strong maximum principle shows that $Q$ is a constant.

We may normalize $u$ so that it attains its maximum 1 at $x_1\in M$ and its minimum $k<1$ at $x_2\in M$. Since $\partial u(x_1)=\partial u(x_2)=0$, we have
\[\frac{K}{2}=\frac{\lambda_1}{4n}u^2(x_1)=Q(x_1)=Q(x_2)=\frac{\lambda_1}{4n}u^2(x_2)=\frac{Kk^2}{2},\]
i.e., $k=-1$, and
\[\frac{K}{2}=Q(x)=|\partial u|^2(x)+\frac{K}{2}u^2(x)\]
for any $x\in M$. In particular, whenever $u^2\neq1$,
\begin{equation}\label{3.22}
\frac{|\partial u|}{\sqrt{1-u^2}}=\sqrt{\frac{K}{2}}.
\end{equation}

Let $\gamma(t)$ be a normal geodesic of shortest length connecting $x_1$ and $x_2$, we may assume $u^2(\gamma)\neq1$ other than the points $x_1$ and $x_2$ without loss of generality. Integrating \eqref{3.22} along $\gamma$, we have
\[d(x_1,x_2)\sqrt{\frac{K}{2}}=\int_\gamma\frac{|\partial u|}{\sqrt{1-u^2}}\geq\frac{1}{\sqrt{2}}\int_{-1}^1\frac{du}{\sqrt{1-u^2}}=\frac{\pi}{\sqrt{2}}.\]

Hence, $D\geq\frac{\pi}{\sqrt{K}}$. Combining with \eqref{3.30}, we have \eqref{1.2xx} that $D=\frac{\pi}{\sqrt{K}}$.

If in addition $\omega$ is K\"ahler, \eqref{1.1x} shows that the real Ricci curvature of the underlying $2n$-dimensional compact Riemannian manifold $M$ satisfies
\[Ric^{\mathbb{R}}\geq(2n-1)Kg.\]
By Cheng's maximum diameter rigidity theorem (see \cite[Theorem 1.3]{Cheng75}), $M$ is isometric to a standard sphere $\mathbb{S}^{2n}$ with constant real sectional curvature $K$. Furthermore, $n=1$ and $(M,\omega)$ is isometric to the complex projective line $\mathbb{CP}^1$ with the Fubini-Study metric (up to scaling).$\hfill\Box$

\section{Li--Yau type estimates}
\label{sec4}

In this section, we prove Theorem~\ref{thm1.2} and Corollary~\ref{cor1.2}, which present Li--Yau type first eigenvalue estimates for the Laplace--de~Rham operator on complete balanced Hermitian manifolds under lower bounds for the holomorphic Ricci curvature of the Strominger--Bismut connection

On a compact balanced Hermitian manifold $(M,\omega)$, we may normalize $u$ so that $u_{\min}=-1$ and $u_{\max}\leq1$ without loss of generality. Let $a>1$ be a constant, $v:=\log(a+u)$ and $P:=|\partial v|^2$. Let $$V=(\partial v)^\sharp=h^{i\bar{k}}\frac{\partial v}{\partial \bar{z}^k}\frac{\partial}{\partial z^i}\in \Gamma(M,T^{1,0}M)$$ be the dual vector field of $\partial v\in\Gamma(M,T^{*1,0}M)$ with respect to the balanced metric $\omega$. Denote by
\begin{equation}\label{4.20}
{^{SB}}\nabla^{1,0}\partial v=t'_{ij}dz^i\otimes dz^j\quad\text{with}\quad \overline{t'_{ij}}=t'_{\bar{i}\bar{j}},
\end{equation}
and
\begin{equation}\label{4.21}
{^{SB}}\nabla^{0,1}\partial v=s'_{\bar{i}j}d\bar{z}^i\otimes dz^j\quad\text{with}\quad s'_{i\bar{j}}=\overline{s'_{\bar{i}j}}.
\end{equation}

\begin{proposition}\label{lem4.1}
On a compact balanced Hermitian manifold $(M,\omega)$, one has
\begin{equation}\label{4.1}
\mathrm{tr}_\omega(\sqrt{-1}\, \partial\bar{\partial}v)=h^{i\bar{j}}s'_{i\bar{j}}=-P-\frac{\lambda_1}{2}+\frac{a\lambda_1}{2(a+u)}.
\end{equation}
\end{proposition}

\begin{proof}
By direct computation, we obtain
\begin{eqnarray}\label{4.2}
\mathrm{tr}_\omega(\sqrt{-1}\, \partial\bar{\partial}v)&=&h^{i\bar{j}}\frac{\partial}{\partial z^i}\big(\frac{1}{a+u}\frac{\partial u}{\partial z^i}\big)\notag\\
&=&-\frac{1}{(a+u)^2}h^{i\bar{j}}\frac{\partial u}{\partial z^i}\frac{\partial u}{\partial \bar{z}^j}+\frac{1}{a+u}h^{i\bar{j}}\frac{\partial^2u}{\partial z^i\partial\bar{z}^j}\notag\\
&=&-P-\frac{\lambda_1}{2}+\frac{a\lambda_1}{2(a+u)},
\end{eqnarray}
where we used Lemma \ref{lem3.1}.

Moreover, we have
\begin{eqnarray}\label{4.19}
s'_{i\bar{j}}&=&\frac{\partial v}{\partial z^i\partial\bar{z}^j}-\overline{{^{SB}}\Gamma_{\bar{i}j}^k}\frac{\partial v}{\partial\bar{z}^k}\notag\\
&=&\frac{1}{a+u}\big(\frac{\partial^2u}{\partial z^i\partial\bar{z}^j}-\overline{{^{SB}}\Gamma_{\bar{i}j}^k}\frac{\partial u}{\partial\bar{z}^k}\big)-\frac{\partial v}{\partial z^i}\frac{\partial v}{\partial\bar{z}^j}\notag\\
&=&\frac{s_{i\bar{j}}}{a+u}-\frac{\partial v}{\partial z^i}\frac{\partial v}{\partial\bar{z}^j}.
\end{eqnarray}

Applying \eqref{3.28} to \eqref{4.19}, we get
\begin{equation}\label{4.2x}
h^{i\bar{j}}s'_{i\bar{j}}=-P-\frac{\lambda_1}{2}+\frac{a\lambda_1}{2(a+u)}.
\end{equation}

\eqref{4.1} follows by combining \eqref{4.2} and \eqref{4.2x}.
\end{proof}

Similar to Proposition \ref{prp3.1}, we can obtain an equality on $P$.
\begin{proposition}\label{prp4.2}
On a compact balanced Hermitian manifold $(M,\omega)$, one has the Bochner type formula that
\begin{eqnarray}\label{4.3}
\mathrm{tr}_{\omega}(\sqrt{-1}\, \partial\bar{\partial}P)&=&\mathcal{R}ic^{SB,\mathbb{C}}(V,\overline{V})-2\mathrm{Re}\langle\partial P,\partial v\rangle\notag\\
&&-\frac{a\lambda_1}{a+u}P+|{^{SB}}\nabla^{1,0}\partial v|^2+|{^{SB}}\nabla^{0,1}\partial v|^2.
\end{eqnarray}
\end{proposition}

\begin{proof}
By the same computation as in \eqref{3.13} and \eqref{3.14}, we can get
\begin{eqnarray}\label{4.4}
\Delta_{\bar{\partial}}P&=&-\sqrt{-1}\, Ric^{SB(2)}(V,\overline{V})+|{^{SB}}\nabla^{1,0}\partial v|^2-|{^{SB}}\nabla^{0,1}\partial v|^2\notag\\
&&+\sqrt{-1}\, \Lambda(\bar{\partial}\{{^{SB}}\nabla^{1,0}\partial v,\partial v\}-\partial\{\partial v,{^{SB}}\nabla^{1,0}\partial v\}).
\end{eqnarray}

Moreover, the same computation as in \eqref{3.15} to \eqref{3.18} implies
\begin{eqnarray}\label{4.5}
&&\sqrt{-1}\, \Lambda\bar{\partial}\{{^{SB}}\nabla^{1,0}\partial v,\partial v\}\notag\\
&=&\sqrt{-1}\, Ric^{SB(1)}(V,\overline{V})+\frac{1}{2}\, T\circ\overline{T}(V,\overline{V})\notag\\
&&-h^{i\bar{k}}\frac{\partial}{\partial z^i}(\mathrm{tr}_\omega(\sqrt{-1}\, \partial\bar{\partial}v))\frac{\partial v}{\partial \bar{z}^k}-|{^{SB}}\nabla^{1,0}\partial v|^2.
\end{eqnarray}
Applying \eqref{4.1} to \eqref{4.5}, we obtain that
\begin{eqnarray}\label{4.6}
\sqrt{-1}\, \Lambda\bar{\partial}\{{^{SB}}\nabla^{1,0}\partial v,\partial v\}&=&\sqrt{-1}\, Ric^{SB(1)}(V,\overline{V})+\frac{1}{2}\, T\circ\overline{T}(V,\overline{V})\notag\\
&&+\langle\partial P,\partial v\rangle+\frac{a\lambda_1P}{2(a+u)}-|{^{SB}}\nabla^{1,0}\partial v|^2.
\end{eqnarray}

By \eqref{4.6}, we have
\begin{eqnarray}\label{4.7}
-\sqrt{-1}\, \Lambda\partial\{\partial v,{^{SB}}\nabla^{1,0}\partial v\}&=&\overline{\sqrt{-1}\, \Lambda\bar{\partial}\{{^{SB}}\nabla^{1,0}\partial v,\partial v\}}\notag\\
&=&\sqrt{-1}\, Ric^{SB(1)}(V,\overline{V})+\frac{1}{2}\, T\circ\overline{T}(V,\overline{V})\notag\\
&&+\overline{\langle\partial P,\partial v\rangle}+\frac{a\lambda_1P}{2(a+u)}-|{^{SB}}\nabla^{1,0}\partial v|^2.
\end{eqnarray}

By applying \eqref{4.6} and \eqref{4.7} to \eqref{4.4} and then using \eqref{6.6} and Lemma \ref{lem3.1}, we get \eqref{4.3}.
\end{proof}

\begin{lemma}\label{lem4.3}
On a compact balanced Hermitian manifold $(M,\omega)$, let $x_3\in M$ be a maximum point of $P$. If $P(x_3)>0$, then
\begin{eqnarray}\label{4.8}
0&\geq&\mathcal{R}ic^{SB,\mathbb{C}}(V,\overline{V})(x_3)-\frac{a\lambda_1}{a+u(x_3)}P(x_3)\notag\\
&&+\frac{2}{2n-1}\left(P(x_3)+\frac{\lambda_1}{2}-\frac{a\lambda_1}{2(a+u(x_3))}\right)^2.
\end{eqnarray}
\end{lemma}

\begin{proof}
Choose a unitary frame $\{\mu_1,\mu_2,\cdots,\mu_n\}\subset T^{1,0}M$ at $x_3$ such that $\mu_1=\frac{V}{|V|}$. Since $x_3$ is a maximum point of $P$, we have $\partial P(x_3)=0$. Moreover,
\[
\mu_1(P)=|\partial v|(t'_{11}+s'_{1\bar{1}}).
\]
Therefore,
\begin{equation}\label{4.30}
s'_{1\bar{1}}(x_3)=-t'_{11}(x_3).
\end{equation}

By \eqref{4.30} and Proposition \ref{lem4.1}, we have
\begin{equation}\label{4.31}
\sum_{\alpha=2}^ns'_{\alpha\bar{\alpha}}(x_3)=t'_{11}(x_3)-P(x_3)-\frac{\lambda_1}{2}+\frac{a\lambda_1}{2(a+u(x_3))}.
\end{equation}
It follows that
\begin{eqnarray}\label{4.32}
&&|{^{SB}}\nabla^{1,0}\partial v|^2+|{^{SB}}\nabla^{0,1}\partial v|^2\notag\\
&\geq&|t'_{11}|^2+|s'_{1\bar{1}}|^2+\frac{1}{n-1}|\sum_{\alpha=2}^ns'_{\alpha\bar{\alpha}}|^2\notag\\
&=&2|t'_{11}|^2+\frac{1}{n-1}\big|t'_{11}-P-\frac{\lambda_1}{2}+\frac{a\lambda_1}{2(a+u)}\big|^2\notag\\
&\geq&\frac{2}{2n-1}\big|P+\frac{\lambda_1}{2}-\frac{a\lambda_1}{2(a+u)}\big|^2
\end{eqnarray}
at $x_3$. 

Note that $$\partial P(x_3)=0\quad\text{and}\quad\mathrm{tr}_{\omega}(\sqrt{-1}\, \partial\bar{\partial}P)(x_3)\leq0.$$
We conclude \eqref{4.8} from Proposition \ref{prp4.2} and \eqref{4.32}.
\end{proof}

We are now ready to prove Theorem \ref{thm1.2} and Corollary \ref{cor1.2}.
\\\\\textbf{Proof of Theorem \ref{thm1.2}.} Since $V=\frac{U}{u+a}$, it follows from \eqref{1.12x} that
\begin{equation}\label{4.17}
\mathcal{R}ic^{SB,\mathbb{C}}(V,\bar{V})\geq-KP.
\end{equation}

Since \(v\) is nonconstant, we have \(P(x_3)>0\). Using \eqref{4.8}, \label{4.17} and the fact that
\[
\big(P+\frac{\lambda}{2}-\frac{a\lambda_1}{2(a+u)}\big)^2\geq P^2+\big(\lambda_1-\frac{a\lambda_1}{a+u}\big)P,
\]
we obtain
\begin{eqnarray}\label{4.15}
P(x)&\leq&P(x_3)\notag\\
&\leq&\frac{2n-1}{2}K+\frac{(2n+1)a\lambda_1}{2(a+u(x_3))}-\lambda_1\notag\\
&\leq&\frac{2n+1}{2}\big(K+\frac{a}{a-1}\lambda_1\big)
\end{eqnarray}
for all $x\in M$.

Since
\[0=(\Delta_du,1)=\lambda_1\int_Mu\frac{\omega^n}{n!},\]
$u$ changes sign and we may normalize $u$ so that $u_{\max}=u(x_4)>0$ and $u_{\min}=u(x_5)=-1$ for some $x_4$ and $x_5\in M$.

Integrating $\sqrt{P}=|\partial\log(a+u)|$ along a minimal geodesic $\gamma$ joining $x_4$ and $x_5$, we have for all $a>1$,
\begin{eqnarray*}
\log\big(\frac{a}{a-1}\big)&\leq&\log\big(\frac{a+u_{\max}}{a-1}\big)\\
&=&\int_\gamma d\log(a+u)\\
&\leq&\int_\gamma\sqrt{2}|\partial\log(a+u)|\\
&\leq&(2n+1)^\frac{1}{2}\big(K+\frac{a}{a-1}\lambda_1\big)^\frac{1}{2}D,
\end{eqnarray*}
where we used \eqref{4.15}, i.e.,
\begin{equation}\label{4.16}
\lambda_1\geq\frac{a-1}{a}\big(\frac{1}{(2n+1)D^2}\big(\log\frac{a}{a-1}\big)^2-K\big).
\end{equation}

Maximizing the right-hand-side of \eqref{4.16} as a function of $a$ by setting
\[a=\frac{1}{1-e^{-\alpha}},\]
with
\[\alpha=1+\big(1+(2n+1)KD^2\big)^\frac{1}{2},\]
we obtain that
\begin{eqnarray*}
\lambda_1&\geq&e^{-\alpha}\big(\frac{\alpha^2}{(2n+1)D^2}-K\big)\notag\\
&\geq&\frac{C_1}{D^2}\exp(-C_2\sqrt{K}D)
\end{eqnarray*}
for some constants $C_1$ and $C_2$ depending only on $n$.$\hfill\Box$
\\\\\textbf{Proof of Corollary \ref{cor1.2}.} Since $K>0$, \cite[Theorem 1.3]{Yang25a+} ensures that \eqref{1.12} implies that
$M$ is compact. Similar to \eqref{4.15}, it follows from \eqref{1.12} and \eqref{4.8} that
\[P(x)\leq\frac{2n+1}{2}\frac{a}{a-1}\lambda_1.\]
By the same arguments as in the proof of Theorem \ref{thm1.2}, taking $K=0$ therein, we obtain \eqref{1.13}.$\hfill\Box$
\section{Zhong--Yang type estimates}
\label{sec5}
In this section, we prove Theorem~\ref{thm1.3} and Corollary~\ref{cor1.3}, yielding Zhong--Yang type first eigenvalue estimates on complete balanced Hermitian manifolds.

On a compact balanced Hermitian manifold $(M,\omega)$, normalize the eigenfunction $u$ that $u_{\max}=1$ and $u_{\min}=-k$ with $0<k\leq 1$ and set $$y=\frac{u-\frac{1-k}{2}}{\frac{1+k}{2}},\ \ \text{and}\ \ b=\frac{1-k}{1+k}\in[0,1).$$
Then we have
\[\Delta_dy=\lambda_1(y+b)\]
with $y_{\max}=1$ and $y_{\min}=-1$.

Let $$Y=(\partial y)^\sharp=h^{i\bar{k}}\frac{\partial y}{\partial \bar{z}^k}\frac{\partial}{\partial z^i}\in \Gamma(M,T^{1,0}M)$$ be the dual vector field of $\partial y\in\Gamma(M,T^{*1,0}M)$ with respect to the balanced metric $\omega$.

Set
\[
M^\circ:=\{x\in M:-1<y(x)<1\}.
\]
Define $\theta:M^\circ\rightarrow \big(-\frac{\pi}{2},\frac{\pi}{2}\big)$ by $\theta=\arcsin y$. Denote by
\[p(x):=\frac{|\partial y|^2}{1-y^2}=|\partial\theta|^2,\]
and for \(t\in(-\pi/2,\pi/2)\), define
\[
F(t):=\max_{\{x\in M^\circ:\theta(x)=t\}}p(x).
\]
as Zhong and Yang did in \cite{ZY84}.

Since $$\partial\theta=\frac{\partial y}{\cos\theta},$$ it follows that
\begin{eqnarray}\label{5.2}
\mathrm{tr}_\omega(\sqrt{-1}\, \partial\bar{\partial}\theta)&=&\frac{1}{\cos\theta}h^{i\bar{j}}\frac{\partial y}{\partial z^i\partial\bar{z}^j}+\frac{\sin\theta}{\cos^2\theta}h^{i\bar{j}}\frac{\partial y}{\partial z^i}\frac{\partial \theta}{\partial\bar{z}^j}\notag\\
&=&-\frac{\Delta_dy}{2\cos\theta}+\frac{\sin\theta}{\cos\theta}h^{i\bar{j}}\frac{\partial \theta}{\partial z^i}\frac{\partial \theta}{\partial\bar{z}^j}\notag\\
&=&-\frac{\lambda_1(\sin\theta+b)}{2\cos\theta}+\frac{\sin\theta}{\cos\theta}|\partial\theta|^2,
\end{eqnarray}
where we used Lemma \ref{lem3.1} in the second equality.

\begin{lemma}\label{lem5.1}
Let $(M,\omega)$ be a compact balanced Hermitian manifold. If
\begin{equation}\label{5.3x}
\mathcal{R}ic^{SB,\mathbb{C}}(Y,\overline{Y}) \ge 0,
\end{equation}
then
\begin{equation}\label{5.3}
F(t)\leq\frac{1+b}{2}\lambda_1,\qquad t\in(-\frac{\pi}{2},\frac{\pi}{2}).
\end{equation}
\end{lemma}

\begin{proof}
By the standard \(\varepsilon\)-regularization
\[
y_\varepsilon:=\frac{y}{1+\varepsilon},
\qquad
b_\varepsilon:=\frac{b}{1+\varepsilon},
\]
and then letting \(\varepsilon\downarrow0\), we may assume that \(p\)
attains its maximum at some point \(x_6\in M^\circ\). Set
\[
\theta_0:=\theta(x_6).
\]
Then
\[
p(x_6)=F(\theta_0)
=\max_{t\in(-\pi/2,\pi/2)}F(t).
\]
Define on \(M\)
\[
\phi(x):=|\partial y|^2(x)-F(\theta_0)\bigl(1-y(x)^2\bigr).
\]
For \(x\in M^\circ\), we have
\[
\phi(x)=\bigl(p(x)-F(\theta_0)\bigr)\bigl(1-y(x)^2\bigr)\leq0,
\]
while \(\phi=0\) on \(M\setminus M^\circ\). Hence \(\phi\) attains
its maximum at \(x_6\). Moreover, at \(x_6\) we have
\begin{equation}\label{5.4}
0=\partial \phi=\{{^{SB}}\nabla^{1,0}\partial y,\partial y\}+\{\partial y,{^{SB}}\nabla^{0,1}\partial y\}+F(\theta_0)\sin2\theta_0\cdot\partial\theta,
\end{equation}
and
\begin{eqnarray}\label{5.5}
0&\geq&\mathrm{tr}_\omega(\sqrt{-1}\, \partial\bar{\partial}\phi)\notag\\
&=&-\Delta_{\bar{\partial}}|\partial y|^2-F(\theta_0)h^{i\bar{j}}\frac{\partial^2\cos^2\theta}{\partial z^i\partial \bar{z}^j}\notag\\
&=&\mathcal{R}ic^{SB,\mathbb{C}}(Y,\overline{Y})-\lambda_1|\partial y|^2+|{^{SB}}\nabla^{1,0}\partial y|^2+|{^{SB}}\nabla^{0,1}\partial y|^2\notag\\
&&-2F(\theta_0)(\sin^2\theta_0-\cos^2\theta_0)h^{i\bar{j}}\frac{\partial\theta}{\partial z^i}\frac{\partial\theta}{\partial \bar{z}^j}+2F(\theta_0)\sin\theta_0\cos\theta_0h^{i\bar{j}}\frac{\partial\theta}{\partial z^i\partial \bar{z}^j}\notag\\
&=&\mathcal{R}ic^{SB,\mathbb{C}}(Y,\overline{Y})-\lambda_1\cos^2\theta_0|\partial \theta|^2+|{^{SB}}\nabla^{1,0}\partial y|^2+|{^{SB}}\nabla^{0,1}\partial y|^2\notag\\
&&-\lambda_1F(\theta_0)(\sin\theta_0+b)\sin\theta_0+2F(\theta_0)\cos^2\theta_0|\partial\theta|^2\notag\\
&=&\mathcal{R}ic^{SB,\mathbb{C}}(Y,\overline{Y})+|{^{SB}}\nabla^{1,0}\partial y|^2+|{^{SB}}\nabla^{0,1}\partial y|^2\notag\\
&&+2F^2(\theta_0)\cos^2\theta_0-\lambda_1F(\theta_0)(1+b\sin\theta_0)
\end{eqnarray}
where we used Lemma \ref{lem3.1} in the first equality, and Proposition \ref{prp3.1} in the second, and \eqref{5.2} in the third.

By \eqref{5.4} and Cauchy-Schwartz inequality, we obtain
\begin{eqnarray*}
F^2(\theta_0)\sin^22\theta_0\cdot|\partial\theta|^2&=&|\{{^{SB}}\nabla^{1,0}\partial y,\partial y\}+\{\partial y,{^{SB}}\nabla^{0,1}\partial y\}|^2\notag\\
&\leq&2(|{^{SB}}\nabla^{1,0}\partial y|^2+|{^{SB}}\nabla^{0,1}\partial y|^2)|\partial y|^2
\end{eqnarray*}
at $x_6$, namely,
\begin{equation}\label{5.6}
|{^{SB}}\nabla^{1,0}\partial y|^2+|{^{SB}}\nabla^{0,1}\partial y|^2\geq2F^2(\theta_0)\sin^2\theta_0
\end{equation}
at $x_6$.

Applying \eqref{5.3x} and \eqref{5.6} to \eqref{5.5}, we obtain
\[
0\geq
2F^2(\theta_0)
-\lambda_1F(\theta_0)
\bigl(1+b\sin\theta_0\bigr).
\]
If \(F(\theta_0)=0\), the conclusion is immediate. Otherwise, dividing
by \(F(\theta_0)\), we obtain
\[
F(\theta_0)
\leq
\frac{\lambda_1}{2}
\bigl(1+b\sin\theta_0\bigr).
\]

It follows that, for every
\(t\in(-\pi/2,\pi/2)\),
\[
F(t)\leq F(\theta_0)
\leq \frac{\lambda_1}{2}(1+b\sin\theta_0)
\leq \frac{1+b}{2}\lambda_1.
\]
\end{proof}

We now proceed to the proof of Theorem \ref{thm1.3}, from which, together with the relation $Y=\frac{2}{1+k}U$, Corollary \ref{cor1.3} follows immediately.
\\\\\textbf{Proof of Theorem \ref{thm1.3}.} As $K>0$, \cite[Theorem 1.3]{Yang25a+} ensures that \eqref{1.14} implies that
$M$ is compact.

Applying the one-dimensional comparison argument of
\cite[Lemmas 3--5]{ZY84} to the preceding local maximum-principle
calculation, with an upper supporting function for \(F\) in place of the
constant \(F(\theta_0)\), yields
\[
2F(\theta)\leq\lambda_1\bigl(1+b\psi(\theta)\bigr).
\]
Consequently,
\begin{equation}\label{5.7}
\sqrt{\lambda_1}\geq \big(\frac{2F(\theta)}{1+b\psi(\theta)}\big)^\frac{1}{2}\geq\frac{|d\theta|}{\sqrt{1+b\psi(\theta)}},
\end{equation}
where $\psi\in C^0\big[-\frac{\pi}{2},\frac{\pi}{2}\big]\cap C^2\big(-\frac{\pi}{2},\frac{\pi}{2}\big)$ is defined in  \cite[Lemmas 4]{ZY84} that $\psi(\frac{\pi}{2}):=1$, $\psi(-\frac{\pi}{2}):=-1$ and
\[\psi(\theta):=\frac{\frac{4}{\pi}(\theta+\cos\theta\sin\theta)-2\sin\theta}{\cos^2\theta},\quad\theta\in\big(-\frac{\pi}{2},\frac{\pi}{2}\big).\]

As in \cite{ZY84}, by integrating \eqref{5.7} and applying Taylor's formula, we obtain that
\begin{eqnarray}\label{5.8}
\sqrt{\lambda_1}D&\geq&\int_{-\frac{\pi}{2}}^{\frac{\pi}{2}}\frac{d\theta}{\sqrt{1+b\psi(\theta)}}\notag\\
&=&\int_0^{\frac{\pi}{2}}\big(\frac{1}{\sqrt{1+b\psi(\theta)}}+\frac{1}{\sqrt{1-b\psi(\theta)}}\big)d\theta\notag\\
&=&\pi\big(1+\sum_{k=1}^\infty\frac{1\times3\times\cdots\times(4k-1)}{2\times 4\times\cdots\times(4k)}C_kb^{2k}\big),
\end{eqnarray}
where $C_k=\frac{2}{\pi}\int_0^{\frac{\pi}{2}}\psi^{2k}(\theta)d\theta$ are positive constants.

Since $b\geq0$, we conclude \eqref{1.15} from \eqref{5.8}.$\hfill\Box$
\section{An integral formula and the proof of Theorem \ref{thm1.4}}
\label{sec6}

The main result of this section is the following integral identity in terms of the Chern connection ${^C}\nabla$, which conditionally extends \cite[Theorem 1.5]{WY25+} from compact K\"ahler manifolds to compact balanced Hermitian manifolds. We denote by $${^C}T_{ij\bar{l}}=h_{k\bar{l}}{^C}T_{ij}^k=\frac{\partial h_{j\bar{l}}}{\partial z^i}-\frac{\partial h_{i\bar{l}}}{\partial z^j},$$
and
\[{^C}T(U,\cdot,\overline{U})=h^{i\bar{a}}h^{b\bar{l}}{^C}T_{ij\bar{l}}\frac{\partial u}{\partial\bar{z}^a}\frac{\partial u}{\partial z^b}dz^j,\qquad{^C}\overline{T}(\overline{U},\cdot,U)=\overline{{^C}T(U,\cdot,\overline{U})}.\]
\begin{theorem}\label{thm6.1}
Let $(M,\omega)$ be a compact balanced Hermitian manifold of complex dimension $n$, if $[\Lambda,\bar{\partial}\omega]\partial u=0$, then the following identity holds.
\begin{eqnarray}\label{1.3}
\lambda_1\int_M|\partial u|^4\frac{\omega^n}{n!}&=&\int_M\Theta(U,\overline{U},U,\overline{U})\frac{\omega^n}{n!}+\int_M|\partial u|^2|\partial\bar{\partial}u|^2\frac{\omega^n}{n!}\notag\\
&&+\|\{{^C}\nabla^{1,0}\partial u,\partial u\}-\frac{\lambda_1}{2} u\partial u\|^2+\|\{{^C}\nabla^{1,0}\partial u,\partial u\}\|^2\notag\\
&&-\mathrm{Re}({^C}T(U,\cdot,\overline{U}),\{{^C}\nabla^{1,0}\partial u,\partial u\})\notag\\
&&+\mathrm{Re}({^C}\overline{T}(\overline{U},\cdot,U),\{{^C}\nabla^{0,1}\partial u,\partial u\}).
\end{eqnarray}
\end{theorem}

\cite[Lemma 2.2]{WY25+} can be generalized to complete Hermitian manifolds with minor modifications by using the identities (see, e.g., \cite[Lemma A.5]{LY12}) that
\begin{equation}\label{3.9}
-\sqrt{-1}\, [\Lambda,\partial]=\bar{\partial}^*+\bar{\tau}^*,
\end{equation}
and
\begin{equation}\label{3.10}
\sqrt{-1}\, [\Lambda,\bar{\partial}]=\partial^*+\tau^*,
\end{equation}
where $\tau^*=-*\bar{\tau}*$ an operator dual to $\tau=[\Lambda,\partial\omega]$, instead of $-\sqrt{-1}\, [\Lambda,\partial]=\bar{\partial}^*$ and $\sqrt{-1}\, [\Lambda,\bar{\partial}]=\partial^*$ in the K\"ahler case:
\begin{lemma}\label{lem2.2}
Given $f\in C^\infty(M,\mathbb{R})$ on a complete Hermitian manifold $(M,\omega)$, one has
\begin{equation}
(\tau^*+\partial^*)(f\partial u\wedge\bar{\partial}u)=f(\Delta_{\bar{\partial}}u)\bar{\partial}u-f\{\partial u,{^C}\nabla^{1,0}\partial u\}-\langle\partial u,\partial f\rangle\bar{\partial}u,
\end{equation}
and
\begin{equation}
(\bar{\tau}^*+\bar{\partial}^*)(f\partial u\wedge\bar{\partial}u)=-f(\Delta_{\partial}u)\partial u+f\{{^C}\nabla^{1,0}\partial u,\partial u\}+\langle\partial f,\partial u\rangle\partial u.
\end{equation}
\end{lemma}

\begin{proposition}\label{lem2.3}
On a compact Hermitian manifold $(M,\omega)$, one has
\begin{eqnarray}\label{2.1}
&&-(\tau^*(\partial u\wedge\bar{\partial}u),\{{^C}\nabla^{0,1}\partial u,\partial u\})\notag\\
&=&(\langle\partial u,\sqrt{-1}\, \bar{\partial}^*\omega\rangle\bar{\partial} u+{^C}\overline{T}(\overline{U},\cdot,U),\{{^C}\nabla^{0,1}\partial u,\partial u\}),
\end{eqnarray}
and
\begin{eqnarray}\label{2.2}
&&-(\bar{\tau}^*(\partial u\wedge\bar{\partial}u),\{{^C}\nabla^{1,0}\partial u,\partial u\})\notag\\
&=&(\langle\bar{\partial}u,\sqrt{-1}\, \partial^*\omega\rangle\partial u-{^C}T(U,\cdot,\overline{U}),\{{^C}\nabla^{1,0}\partial u,\partial u\}).
\end{eqnarray}
\end{proposition}

\begin{proof}
Since the operator $\tau^*$ is dual to $\tau$,
\begin{eqnarray}\label{2.3}
&&-(\tau^*(\partial u\wedge\bar{\partial}u),\{{^C}\nabla^{0,1}\partial u,\partial u\})\notag\\
&=&(\partial u\wedge\bar{\partial}u,-\tau\{{^C}\nabla^{0,1}\partial u,\partial u\})\notag\\
&=&\int_M\langle\partial u\wedge\bar{\partial}u,-[\Lambda,\partial\omega]\{{^C}\nabla^{0,1}\partial u,\partial u\}\rangle\frac{\omega^n}{n!}.
\end{eqnarray}

Moreover, we have
\begin{eqnarray*}
&&-[\Lambda,\partial\omega]\{{^C}\nabla^{0,1}\partial u,\partial u\}\notag\\
&=&-\Lambda(\partial\omega\wedge\{{^C}\nabla^{0,1}\partial u,\partial u\})\notag\\
&=&\sqrt{-1}\, h^{a\bar{b}}\iota_{\bar{b}}\iota_a\big(\sqrt{-1}\, \frac{\partial h_{p\bar{q}}}{\partial z^s}h^{i\bar{k}}\frac{\partial^2u}{\partial z^i\partial\bar{z}^j}\frac{\partial u}{\partial\bar{z}^k}dz^s\wedge dz^p\wedge d\bar{z}^q\wedge d\bar{z}^j\big)\notag\\
&=&h^{a\bar{b}}\iota_{\bar{b}}\big({^C}T_{pa\bar{q}}h^{i\bar{k}}\frac{\partial^2u}{\partial z^i\partial\bar{z}^j}\frac{\partial u}{\partial\bar{z}^k} dz^p\wedge d\bar{z}^q\wedge d\bar{z}^j\big)\notag\\
&=&-h^{a\bar{b}}{^C}T_{pa\bar{b}}h^{i\bar{k}}\frac{\partial^2u}{\partial z^i\partial\bar{z}^j}\frac{\partial u}{\partial \bar{z}^k} dz^p\wedge d\bar{z}^j\notag\\
&&+h^{a\bar{b}}{^C}T_{pa\bar{q}}h^{i\bar{k}}\frac{\partial^2u}{\partial z^i\partial\bar{z}^b}\frac{\partial u}{\partial\bar{z}^k} dz^p\wedge d\bar{z}^q.
\end{eqnarray*}

By \eqref{3.29}, we have
\[\bar{\partial}^*\omega=\sqrt{-1}\, h^{i\bar{j}}{^C}T_{pi\bar{j}}dz^p.\]
Therefore,
\begin{eqnarray}\label{2.4}
&&\langle\partial u\wedge\bar{\partial}u,-[\Lambda,\partial\omega]\{{^C}\nabla^{0,1}\partial u,\partial u\}\rangle\notag\\
&=&-h^{l\bar{p}}h^{j\bar{q}}h^{b\bar{a}}\overline{{^C}T_{pa\bar{b}}}h^{k\bar{i}}\frac{\partial^2u}{\partial z^j\partial\bar{z}^i}\frac{\partial u}{\partial z^k}\frac{\partial u}{\partial z^l}\frac{\partial u}{\partial\bar{z}^q}\notag\\
&&+h^{l\bar{p}}h^{q\bar{j}}h^{b\bar{a}}\overline{{^C}T_{pa\bar{q}}}h^{k\bar{i}}\frac{\partial^2u}{\partial z^b\partial\bar{z}^i}\frac{\partial u}{\partial z^k}\frac{\partial u}{\partial z^l}\frac{\partial u}{\partial\bar{z}^j}\notag\\
&=&\langle\langle\partial u,\sqrt{-1}\, \bar{\partial}^*\omega\rangle\cdot\bar{\partial}u+{^C}\overline{T}(\overline{U},\cdot,U),\{{^C}\nabla^{0,1}\partial u,\partial u\}\rangle.
\end{eqnarray}
\eqref{2.1} follows by applying \eqref{2.4} into \eqref{2.3}. \eqref{2.2} can be proved in a similar way.
\end{proof}

\begin{lemma}\label{lem6.4}
Given $f\in C^\infty(M,\mathbb{R})$ on a compact balanced Hermitian manifold $(M,\omega)$, one has
\begin{equation}\label{6.11}
\partial^*\partial\bar{\partial}f=\bar{\partial}\partial^*\partial f-\partial^*\bar{\tau}\partial f.
\end{equation}
\end{lemma}

\begin{proof}
For any $\varphi\in\Gamma(M,T^{*0,1}M)$, we have
\begin{eqnarray}\label{6.12}
&&(\bar{\partial}\partial^*\partial f-\partial^*\partial\bar{\partial}f,\varphi)\notag\\
&=&(\partial f,\partial\bar{\partial}^*\varphi)+(\partial f,\bar{\partial}^*\partial\varphi)\notag\\
&=&-(\partial f,\partial(\sqrt{-1}[\Lambda,\partial]+\bar{\tau}^*)\varphi)\notag\\
&&-(\partial f,(\sqrt{-1}[\Lambda,\partial]+\bar{\tau}^*)\partial\varphi)\notag\\
&=&-(\partial f,\sqrt{-1}\partial(\Lambda\partial\varphi)-\sqrt{-1}\partial(\Lambda\partial\varphi))\notag\\
&&-(\partial f,\partial\bar{\tau}^*\varphi+\bar{\tau}^*\partial\varphi)\notag\\
&=&-(\bar{\tau}\partial^*\partial f,\varphi)+(\partial^*\bar{\tau}\partial f,\varphi).
\end{eqnarray}
where we used \eqref{3.9} in the second equality.

It follows from \eqref{3.29} and $\partial^*\omega=0$ that
\begin{equation}\label{6.13}
\bar{\tau}\partial^*\partial f=\overline{[\Lambda,\partial\omega]}\partial^*\partial f=\partial^*\partial f\Lambda(\bar{\partial}\omega)=\sqrt{-1}\partial^*\partial f\wedge\partial^*\omega=0.
\end{equation}

We conclude \eqref{6.11} by applying \eqref{6.13} to \eqref{6.12}.
\end{proof}

\begin{proposition}\label{prp6.2}
On a compact balanced Hermitian manifold, if $[\Lambda,\bar{\partial}\omega]\partial u=0$, then
\begin{eqnarray}\label{3.4}
\lambda_1\int_M|\partial u|^4\frac{\omega^n}{n!}&=&\int_M|\partial u|^2|\partial\bar{\partial}u|^2\frac{\omega^n}{n!}+\|\{\partial u,{^C}\nabla^{0,1}\partial u\}\|^2\notag\\
&&+\mathrm{Re}(\{{^C}\nabla^{1,0}\partial u,\partial u\},\{\partial u,{^C}\nabla^{0,1}\partial u\})+\frac{\lambda_1^2}{4}\|u\partial u\|^2\notag\\
&&-\frac{\lambda_1}{2}\mathrm{Re}(\{{^C}\nabla^{1,0}\partial u,\partial u\}+\{\partial u,{^C}\nabla^{0,1}\partial u\},u\partial u).
\end{eqnarray}
\end{proposition}

\begin{proof}
By Lemma \ref{lem3.1}, we obtain that
\begin{equation}\label{6.9}
(\partial|\partial u|^2,\Delta_{\bar{\partial}}u\cdot\partial u)=\frac{\lambda_1}{2}(\{{^C}\nabla^{1,0}\partial u,\partial u\}+\{\partial u,{^C}\nabla^{0,1}\partial u\},u\partial u).
\end{equation}

On the other hand, we have
\begin{eqnarray}\label{6.10}
&&(\partial|\partial u|^2,\Delta_{\bar{\partial}}u\cdot\partial u)\notag\\
&=&\frac{\lambda_1}{2}(|\partial u|^2,\partial^*(u\partial u))\notag\\
&=&\frac{\lambda_1}{2}(|\partial u|^2,\sqrt{-1}\, \Lambda\bar{\partial}(u\partial u))\notag\\
&=&-\frac{\lambda_1}{2}(|\partial u|^2,\mathrm{tr}_\omega(\sqrt{-1}\, \partial u\wedge\bar{\partial}u))-\frac{\lambda_1}{2}(|\partial u|^2,u\cdot\mathrm{tr}_\omega(\sqrt{-1}\, \partial \bar{\partial}u))\notag\\
&=&-\frac{\lambda_1}{2}\int_M|\partial u|^4\frac{\omega^n}{n!}+\frac{\lambda_1^2}{4}\|u\partial u\|^2,
\end{eqnarray}
where we used Lemma \ref{lem3.1} in the first and last equalities, and \eqref{3.7x} in the second.

By \eqref{6.9} and \eqref{6.10} and taking conjugate, we have
\begin{eqnarray}\label{3.2}
&&\frac{\lambda_1}{2}\int_M|\partial u|^4\frac{\omega^n}{n!}\notag\\
&=&\frac{\lambda_1^2}{4}\|u\partial u\|^2-\frac{\lambda_1}{2}\mathrm{Re}(\{{^C}\nabla^{1,0}\partial u,\partial u\}+\{\partial u,{^C}\nabla^{0,1}\partial u\},u\partial u).
\end{eqnarray}

Since $[\Lambda,\bar{\partial}\omega]\partial u=0$, one has
\begin{equation}\label{6.14}
\partial^*\bar{\tau}\partial u=\partial^*([\Lambda,\bar{\partial}\omega]\partial u)=0.
\end{equation}
It follows from \eqref{6.14}, Lemma \ref{lem6.4} and Lemma \ref{lem3.1} that
\[\partial^*\partial\bar{\partial}u=\bar{\partial}\partial^*\partial u=\frac{1}{2}\bar{\partial}(\Delta_du)=\frac{\lambda_1}{2}\bar{\partial}u.\]
Therefore,
\begin{eqnarray*}
&&\int_M|\partial u|^2|\partial\bar{\partial}u|^2\frac{\omega^n}{n!}\notag\\
&=&(|\partial u|^2\cdot\partial\bar{\partial}u,\partial\bar{\partial}u)\notag\\
&=&(\partial(|\partial u|^2\bar{\partial}u)-\partial|\partial u|^2\wedge\bar{\partial}u,\partial\bar{\partial}u)\notag\\
&=&(|\partial u|^2\bar{\partial}u,\partial^*\partial\bar{\partial}u)-(\partial|\partial u|^2,\{\partial u,{^C}\nabla^{0,1}\partial u\})\notag\\
&=&\frac{\lambda_1}{2}\int_M|\partial u|^4\frac{\omega^n}{n!}-(\{{^C}\nabla^{1,0}\partial u,\partial u\},\{\partial u,{^C}\nabla^{0,1}\partial u\})\notag\\
&&-\|\{\partial u,{^C}\nabla^{0,1}\partial u\}\|^2.
\end{eqnarray*}
Taking conjugate, we have
\begin{eqnarray}\label{3.2x}
\frac{\lambda_1}{2}\int_M|\partial u|^4\frac{\omega^n}{n!}&=&\int_M|\partial u|^2|\partial\bar{\partial}u|^2\frac{\omega^n}{n!}+\|\{\partial u,{^C}\nabla^{0,1}\partial u\}\|^2\notag\\
&&+\mathrm{Re}(\{{^C}\nabla^{1,0}\partial u,\partial u\},\{\partial u,{^C}\nabla^{0,1}\partial u\}).
\end{eqnarray}

Summing up \eqref{3.2} and \eqref{3.2x}, we get \eqref{3.4}.
\end{proof}

Based on the preceding identities, we proceed to the proof of Theorem \ref{thm6.1}.
\\\textbf{Proof of Theorem \ref{thm6.1}.}
Consider the $(1,1)$-form
\begin{eqnarray}
\alpha&=&\partial\{{^C}\nabla^{0,1}\partial u,\partial u\}+\bar{\partial}\{{^C}\nabla^{1,0}\partial u,\partial u\}\notag\\
&=&-\{h^{b\bar{l}}\Theta_{i\bar{j}k\bar{l}}\frac{\partial u}{\partial z^b}dz^i\wedge d\bar{z}^j\otimes dz^k,\frac{\partial u}{\partial z^a}dz^a\}\notag\\
&&-\{{^C}\nabla^{1,0}\partial u,{^C}\nabla^{1,0}\partial u\}-\{{^C}\nabla^{0,1}\partial u,{^C}\nabla^{0,1}\partial u\}.
\end{eqnarray}
It follows that
\begin{eqnarray}\label{1.6}
(\partial u\wedge\bar{\partial} u,\alpha)&=&-\int_M\overline{\Theta(U,\overline{U},U,\overline{U})}\frac{\omega^n}{n!}-\|\{{^C}\nabla^{1,0}\partial u,\partial u\}\|^2\notag\\
&&+\|\{{^C}\nabla^{0,1}\partial u,\partial u\}\|^2.
\end{eqnarray}

On the other hand, using Lemma \ref{lem2.2}, Proposition \ref{lem2.3} and $d^*\omega=0$, we have
\begin{eqnarray}\label{1.7}
&&(\partial u\wedge\bar{\partial} u,\alpha)\notag\\
&=&(\partial u\wedge\bar{\partial} u,\partial\{{^C}\nabla^{0,1}\partial u,\partial u\}+\bar{\partial}\{{^C}\nabla^{1,0}\partial u,\partial u\})\notag\\
&=&(\partial^*(\partial u\wedge\bar{\partial} u),\{{^C}\nabla^{0,1}\partial u,\partial u\})+(\bar{\partial}^*(\partial u\wedge\bar{\partial} u),\{{^C}\nabla^{1,0}\partial u,\partial u\})\notag\\
&=&-(\tau^*(\partial u\wedge\bar{\partial} u),\{{^C}\nabla^{0,1}\partial u,\partial u\})-(\bar{\tau}^*(\partial u\wedge\bar{\partial} u),\{{^C}\nabla^{1,0}\partial u,\partial u\})\notag\\
&&+(\Delta_{\bar{\partial}}u\cdot\bar{\partial}u,\{{^C}\nabla^{0,1}\partial u,\partial u\})-(\{\partial u,{^C}\nabla^{1,0}\partial u\},\{{^C}\nabla^{0,1}\partial u,\partial u\})\notag\\
&&-(\Delta_{\partial}u\cdot\partial u,\{{^C}\nabla^{1,0}\partial u,\partial u\})+\|\{{^C}\nabla^{1,0}\partial u,\partial u\}\|^2\notag\\
&=&({^C}\overline{T}(\overline{U},\cdot,U),\{{^C}\nabla^{0,1}\partial u,\partial u\})-({^C}T(U,\cdot,\overline{U}),\{{^C}\nabla^{1,0}\partial u,\partial u\})\notag\\
&&+\frac{\lambda_1}{2}(\{\partial u,{^C}\nabla^{0,1}\partial u\},u\partial u)-(\{\partial u,{^C}\nabla^{1,0}\partial u\},\{{^C}\nabla^{0,1}\partial u,\partial u\})\notag\\
&&-\frac{\lambda_1}{2}(u\partial u,\{{^C}\nabla^{1,0}\partial u,\partial u\})+\|\{{^C}\nabla^{1,0}\partial u,\partial u\}\|^2.
\end{eqnarray}

Combining \eqref{1.6} and \eqref{1.7} and taking conjugate, we obtain
\begin{eqnarray}\label{3.1}
0&=&\int_M\Theta(U,\overline{U},U,\overline{U})\frac{\omega^n}{n!}+2\|\{{^C}\nabla^{1,0}\partial u,\partial u\}\|^2\notag\\
&&-\|\{{^C}\nabla^{0,1}\partial u,\partial u\}\|^2-\mathrm{Re}(\{\partial u,{^C}\nabla^{1,0}\partial u\},\{{^C}\nabla^{0,1}\partial u,\partial u\})\notag\\
&&+\frac{\lambda_1}{2}\mathrm{Re}(\{\partial u,{^C}\nabla^{0,1}\partial u\},u\partial u)-\frac{\lambda_1}{2}\mathrm{Re}(u\partial u,\{{^C}\nabla^{1,0}\partial u,\partial u\})\notag\\
&&-\mathrm{Re}({^C}T(U,\cdot,\overline{U}),\{{^C}\nabla^{1,0}\partial u,\partial u\})\notag\\
&&+\mathrm{Re}({^C}\overline{T}(\overline{U},\cdot,U),\{{^C}\nabla^{0,1}\partial u,\partial u\}).
\end{eqnarray}

The summation of \eqref{3.4} and \eqref{3.1} gives \eqref{1.3}.$\hfill\Box$

As an application of Theorem \ref{thm6.1}, we have the following first eigenvalue estimate.
\begin{proposition}\label{prp6.1}
Let $(M,\omega)$ be a compact balanced Hermitian manifold of complex dimension $n$, if $[\Lambda,\bar{\partial}\omega]\partial u=0$, then
\begin{equation}\label{1.4}
\lambda_1\int_M|\partial u|^4\frac{\omega^n}{n!}\geq\int_MR^{SB,\mathbb{C}}(U,\overline{U},U,\overline{U})\frac{\omega^n}{n!}+\frac{1}{2}\|{^C}T(U,\cdot,\overline{U})\|^2.
\end{equation}
\end{proposition}

\begin{proof}
By \eqref{6.0}, we obtain that
\begin{equation}\label{1.8}
R^{SB,\mathbb{C}}(U,\overline{U},U,\overline{U})=\Theta(U,\overline{U},U,\overline{U})-|{^C}T(U,\cdot,\overline{U})|^2.
\end{equation}

Note that
\begin{eqnarray}\label{1.10}
&&\|\{{^C}\nabla^{1,0}\partial u,\partial u\}\|^2-\mathrm{Re}({^C}T(U,\cdot,\overline{U}),\{{^C}\nabla^{1,0}\partial u,\partial u\})+\frac{1}{4}\|{^C}T(U,\cdot,\overline{U})\|^2\notag\\
&=&\|\frac{1}{2}{^C}T(U,\cdot,\overline{U})-\{{^C}\nabla^{1,0}\partial u,\partial u\}\|^2\geq0,
\end{eqnarray}
and
\begin{eqnarray}\label{1.10x}
&&|\partial u|^2|\partial\bar{\partial}u|^2+\mathrm{Re}({^C}\overline{T}(\overline{U},\cdot,U),\{{^C}\nabla^{0,1}\partial u,\partial u\})+\frac{1}{4}\|{^C}T(U,\cdot,\overline{U})\|^2\notag\\
&\geq&\|\{{^C}\nabla^{0,1}\partial u,\partial u\}\|^2+\mathrm{Re}({^C}\overline{T}(\overline{U},\cdot,U),\{{^C}\nabla^{0,1}\partial u,\partial u\})+\frac{1}{4}\|{^C}\overline{T}(\overline{U},\cdot,U)\|^2\notag\\
&=&\|\frac{1}{2}{^C}\overline{T}(\overline{U},\cdot,U)+\{{^C}\nabla^{0,1}\partial u,\partial u\}\|^2\geq0.
\end{eqnarray}

Applying \eqref{1.8}, \eqref{1.10} and \eqref{1.10x} to \eqref{1.3}, we get \eqref{1.4}.
\end{proof}

We complete the proof of Theorem \ref{thm1.4} by applying Proposition \ref{prp6.1}.
\\\textbf{Proof of Theorem \ref{thm1.4}.}  As $K>0$, \cite[Theorem 1.4]{Yang25a+} ensures that \eqref{6.8} implies that
$M$ is compact.

Let $X=U+\overline{U}\in\Gamma(M,T_{\mathbb{R}}M)$. At every point where $U\neq0$, set $$JX=\sqrt{-1}\, (U-\overline{U}).$$
It is clear that
\[g(X,X)=g(U+\overline{U},U+\overline{U})=2h(U,\overline{U}),\]
and
\begin{eqnarray*}
&&R^{SB,\mathbb{R}}(JX,X,X,JX)\notag\\
&=&R^{SB,\mathbb{C}}(\sqrt{-1}\, U-\sqrt{-1}\, \overline{U},U+\overline{U},U+\overline{U}\sqrt{-1}\, U-\sqrt{-1}\, \overline{U})\notag\\
&=&-R^{SB,\mathbb{C}}(U,U,U,U)+R^{SB,\mathbb{C}}(U,U,U,\overline{U})-R^{SB,\mathbb{C}}(U,U,\overline{U},U)\notag\\
&&+R^{SB,\mathbb{C}}(U,U,\overline{U},\overline{U})-R^{SB,\mathbb{C}}(U,\overline{U},U,U)+R^{SB,\mathbb{C}}(U,\overline{U},U,\overline{U})\notag\\
&&-R^{SB,\mathbb{C}}(U,\overline{U},\overline{U},U)+R^{SB,\mathbb{C}}(U,\overline{U},\overline{U},\overline{U})+R^{SB,\mathbb{C}}(\overline{U},U,U,U)\notag\\
&&-R^{SB,\mathbb{C}}(\overline{U},U,U,\overline{U})+R^{SB,\mathbb{C}}(\overline{U},U,\overline{U},U)-R^{SB,\mathbb{C}}(\overline{U},U,\overline{U},\overline{U})\notag\\
&&+R^{SB,\mathbb{C}}(\overline{U},\overline{U},U,U)-R^{SB,\mathbb{C}}(\overline{U},\overline{U},U,\overline{U})+R^{SB,\mathbb{C}}(\overline{U},\overline{U},\overline{U},U)\notag\\
&&-R^{SB,\mathbb{C}}(\overline{U},\overline{U},\overline{U},\overline{U})\notag\\
&=&4R^{SB,\mathbb{C}}(U,\overline{U},U,\overline{U}).
\end{eqnarray*}

It follows that
\begin{equation}\label{6.7}
K\leq\mathrm{HSC}^{SB}(X)=\frac{R^{SB,\mathbb{R}}(JX,X,X,JX)}{g(X,X)^2}=\frac{R^{SB,\mathbb{C}}(U,\overline{U},U,\overline{U})}{h(U,\overline{U})^2}.
\end{equation}

Applying \eqref{6.7} to Proposition \ref{prp6.1}, we get \eqref{6.15} that $\lambda_1\geq K$.$\hfill\Box$
\\\\\textbf{Proof of Corollary \ref{cor1.5}.} By direct computation, one has
\begin{eqnarray}
[\Lambda,\bar{\partial}\omega]\partial u&=&\bar{\tau}\partial u\notag\\
&=&\Lambda(\bar{\partial}\omega\wedge\partial u)\notag\\
&=&\Lambda(\bar{\partial}\omega)\wedge\partial u-\overline{{^C}T_{iq\bar{p}}U^i}dz^p\wedge d\bar{z}^q\notag\\
&=&\sqrt{-1}\partial^*\omega\wedge\partial u+\overline{h_{k\bar{p}}{^{SB}}T_{iq}^kU^i}dz^p\wedge d\bar{z}^q\notag\\
&=&0,
\end{eqnarray}
and hence Propositions \ref{prp6.2} and \ref{prp6.1} hold. The proof of this corollary is contained in the proof of Theorem \ref{thm1.4}.$\hfill\Box$
\section*{Acknowledgements}
The author wishes to express his gratitude to Professor Kefeng Liu for unwavering guidance and generous encouragement, and to Professors Yunhui Wu, Xiaokui Yang and Tao Zheng for illuminating discussions and kind support.

\section*{Statements and Declarations}

\noindent\textbf{Funding.}
This work was supported by the Scientific Research Foundation of
Chongqing University of Technology (Grant No.\ 2026ZDZ012) and the
Youth Project of the Science and Technology Research Program of
Chongqing Education Commission of China
(Grant No.\ KJQN202601108).

\noindent\textbf{Competing interests.}
The author has no competing interests to declare.

\noindent\textbf{Data availability.}
Data sharing is not applicable to this article as no datasets were generated
or analyzed during the current study.

\noindent\textbf{Use of AI-assisted tools.}
During the preparation of this manuscript, the author used ChatGPT
(OpenAI) for language editing and manuscript-review assistance.
The author independently verified all mathematical arguments, derivations,
and references and takes full responsibility for the final content.

\end{document}